\documentclass[11pt,oneside]{amsart}

\usepackage{fullpage}
\usepackage{amsthm,amsfonts,amssymb,amsmath,amsxtra}
\usepackage[all]{xy}
\SelectTips{cm}{}
\usepackage{tikz}
\usepackage{verbatim}
\usepackage{mathrsfs}

\RequirePackage{xspace}
\RequirePackage{etoolbox}
\RequirePackage{varwidth}
\RequirePackage{enumitem}
\RequirePackage{tensor}
\RequirePackage{mathtools}
\RequirePackage{longtable}
\RequirePackage{multirow}
\RequirePackage{makecell}
\RequirePackage{bigints}

\usepackage[backref=page]{hyperref} %

\usepackage[T2A,T1]{fontenc}
\DeclareSymbolFont{cyrillic}{T2A}{cmr}{m}{n}
\DeclareMathSymbol{\Sha}{\mathalpha}{cyrillic}{216}

\newcommand{\sH}{\ensuremath{\mathscr{H}}\xspace}

\newcommand{\sS}{\ensuremath{\mathscr{S}}\xspace}

\newcommand{\sV}{\ensuremath{\mathscr{V}}\xspace}

\DeclareMathOperator{\Aut}{Aut}

\newcommand{\Ch}{{\mathrm{CH}}}
\newcommand{\aCh}{{\widehat{\mathrm{CH}}}}

\DeclareMathOperator{\diag}{diag}
\newcommand{\disc}{{\mathrm{disc}}}

\newcommand{\GL}{\mathrm{GL}}

\newcommand{\GSpin}{\mathrm{GSpin}}

\DeclareMathOperator{\Hom}{Hom}

\let\Im\relax

\DeclareMathOperator{\Im}{Im}
\newcommand{\Ind}{{\mathrm{Ind}}}

\DeclareMathOperator{\Int}{\ensuremath{\mathrm{Int}}\xspace}

\newcommand{\Mp}{{\mathrm{Mp}}}

\DeclareMathOperator{\ord}{ord}

\DeclareMathOperator{\rank}{rank}

\renewcommand{\Re}{{\mathrm{Re}}}

\newcommand{\Sh}{\mathrm{Sh}}

\newcommand{\SL}{{\mathrm{SL}}}

\DeclareMathOperator{\Spec}{Spec}
\DeclareMathOperator{\Spf}{Spf}

\newcommand{\SO}{{\mathrm{SO}}}
\newcommand{\Sp}{{\mathrm{Sp}}}

\newcommand{\val}{{\mathrm{val}}}

\DeclareMathOperator{\Sym}{Sym}

\DeclareMathOperator{\tr}{tr}

\renewcommand{\U}{\mathrm{U}}
\renewcommand{\O}{\mathrm{O}}

\DeclareMathOperator{\vol}{vol}

\newcommand{\lra}{\longrightarrow}

\newenvironment{altenumerate}
   {\begin{list}
      {(\theenumi) }
      {\usecounter{enumi}
       \setlength{\labelwidth}{0pt}
       \setlength{\labelsep}{0pt}
       \setlength{\leftmargin}{0pt}
       \setlength{\itemsep}{\the\smallskipamount}
       \renewcommand{\theenumi}{\roman{enumi}}
      }}
   {\end{list}}

\setitemize[0]{leftmargin=*,itemsep=\the\smallskipamount}
\setenumerate[0]{leftmargin=*,itemsep=\the\smallskipamount}

\renewcommand{\to}{%
   \ifbool{@display}{\longrightarrow}{\rightarrow}%
   }
\let\shortmapsto\mapsto
\renewcommand{\mapsto}{%
   \ifbool{@display}{\longmapsto}{\shortmapsto}%
   }
\newlength{\olen}
\newlength{\ulen}
\newlength{\xlen}
\newcommand{\xra}[2][]{%
   \ifbool{@display}%
      {\settowidth{\olen}{$\overset{#2}{\longrightarrow}$}%
       \settowidth{\ulen}{$\underset{#1}{\longrightarrow}$}%
       \settowidth{\xlen}{$\xrightarrow[#1]{#2}$}%
       \ifdimgreater{\olen}{\xlen}%
          {\underset{#1}{\overset{#2}{\longrightarrow}}}%
          {\ifdimgreater{\ulen}{\xlen}%
             {\underset{#1}{\overset{#2}{\longrightarrow}}}
             {\xrightarrow[#1]{#2}}}}%
      {\xrightarrow[#1]{#2}}
   }
\makeatother
\newcommand{\xyra}[2][]{%
   \settowidth{\xlen}{$\xrightarrow[#1]{#2}$}%
   \ifbool{@display}%
      {\settowidth{\olen}{$\overset{#2}{\longrightarrow}$}%
       \settowidth{\ulen}{$\underset{#1}{\longrightarrow}$}%
       \ifdimgreater{\olen}{\xlen}%
          {\mathrel{\xymatrix@M=.12ex@C=3.2ex{\ar[r]^-{#2}_-{#1} &}}}%
          {\ifdimgreater{\ulen}{\xlen}%
             {\mathrel{\xymatrix@M=.12ex@C=3.2ex{\ar[r]^-{#2}_-{#1} &}}}
             {\mathrel{\xymatrix@M=.12ex@C=\the\xlen{\ar[r]^-{#2}_-{#1} &}}}}}%
      {\mathrel{\xymatrix@M=.12ex@C=\the\xlen{\ar[r]^-{#2}_-{#1} &}}}%
   }
\makeatletter
\newcommand{\xla}[2][]{%
   \ifbool{@display}%
      {\settowidth{\olen}{$\overset{#2}{\longleftarrow}$}%
       \settowidth{\ulen}{$\underset{#1}{\longleftarrow}$}%
       \settowidth{\xlen}{$\xleftarrow[#1]{#2}$}%
       \ifdimgreater{\olen}{\xlen}%
          {\underset{#1}{\overset{#2}{\longleftarrow}}}%
          {\ifdimgreater{\ulen}{\xlen}%
             {\underset{#1}{\overset{#2}{\longleftarrow}}}
             {\xleftarrow[#1]{#2}}}}%
      {\xleftarrow[#1]{#2}}
   }
\newcommand{\isoarrow}{%
   \ifbool{@display}{\overset{\sim}{\longrightarrow}}{\xrightarrow\sim}%
   }
\renewcommand{\lra}{%
   \ifbool{@display}{\longleftrightarrow}{\leftrightarrow}%
   }

\newcommand{\OFb}{{O_{\breve F}}}

\newcommand{\Herm}{\mathrm{Herm}}
\newcommand{\rd}{\mathrm{d}}
\newcommand{\Rep}{\mathrm{Rep}}
\newcommand{\iden}{\langle1\rangle}

\newcommand{\Den}{\mathrm{Den}}

\newcommand{\pDen}{\partial\mathrm{Den}}

\newcommand{\sx}{\mathsf{x}}
\newcommand{\sy}{\mathsf{y}}

\DeclareSymbolFont{sfletters}{OML}{cmbrm}{m}{it}
\DeclareMathSymbol{\stau}{\mathord}{sfletters}{"1C}
\newcommand{\sz}{\tau}

\DeclareFontFamily{U}{matha}{\hyphenchar\font45}
\DeclareFontShape{U}{matha}{m}{n}{
      <5> <6> <7> <8> <9> <10> gen * matha
      <10.95> matha10 <12> <14.4> <17.28> <20.74> <24.88> matha12
      }{}
\DeclareSymbolFont{matha}{U}{matha}{m}{n}
\DeclareFontFamily{U}{mathx}{\hyphenchar\font45}
\DeclareFontShape{U}{mathx}{m}{n}{
      <5> <6> <7> <8> <9> <10>
      <10.95> <12> <14.4> <17.28> <20.74> <24.88>
      mathx10
      }{}
\DeclareSymbolFont{mathx}{U}{mathx}{m}{n}

\DeclareMathSymbol{\obot}         {2}{matha}{"6B}

\newtheorem{theorem}[subsubsection]{Theorem}
\newtheorem{proposition}[subsubsection]{Proposition}

\newtheorem {conjecture}[subsubsection]{Conjecture}
\newtheorem{corollary}[subsubsection]{Corollary}

\theoremstyle{definition}
\newtheorem{definition}[subsubsection]{Definition}
\newtheorem{example}[subsubsection]{Example}

\newtheorem{remark}[subsubsection]{Remark}

\numberwithin{equation}{subsubsection}
\newtheorem{assumption}[subsubsection]{Assumption}

\newcommand{\m}{n}
\renewcommand{\i}{m}
\newcommand{\Gen}{\mathrm{Gen}}
\newcommand{\KL}{K}

\newcommand{\DL}{\mathop {\operator@font DL}\nolimits}
\newcommand{\UU}{\mathrm{U}}
\newcommand{\OO}{\mathrm{O}}

\title{From sum of two squares to arithmetic Siegel--Weil formulas}

\author[Chao Li]{Chao Li}
\address{Columbia University, Department of Mathematics, 2990 Broadway,	New York, NY 10027, USA}
\email{chaoli@math.columbia.edu}

\date{January 11, 2023}


\dedicatory{In loving memory of my mother Xiaoping Mao (1965--2022).}

\thanks{The author would like to thank Benedict Gross, Yifeng Liu, Murilo Zanarella, Wei Zhang and the anonymous referee for helpful comments. The author's work is partially supported by the NSF grant DMS-2101157.}

\begin{document}

\maketitle{}

\begin{abstract}
The main goal of this expository article is to survey recent progress on the arithmetic Siegel--Weil formula and its applications. We begin with the classical sum of two squares problem and put it in the context of the Siegel--Weil formula. We then motivate the geometric and arithmetic Siegel--Weil formula using the classical example of the product of modular curves. After explaining the recent result on the arithmetic Siegel--Weil formula for Shimura varieties of arbitrary dimension, we discuss some aspects of the proof and its application to the arithmetic inner product formula and the Beilinson--Bloch conjecture. Rather than intended to be a complete survey of this vast field, this article focuses more on examples and background to provide easier access to several recent works by the author with W. Zhang \cite{LZ,LZ22} and Y. Liu \cite{LL2020, LL2021}.
\end{abstract}

\tableofcontents

\section{Sum of two squares}

\subsection{Which prime $p$ can be written as the sum of two squares?}
For the first few primes we easily find that $$5=1^2+2^2, \quad 13=2^2+3^2,\quad 17=1^2+4^2,\quad 29=2^2+5^2$$ are sum of two squares, while other primes like $3, 7, 11, 19,23$ are not. The answer seems to depend on the residue class of $p$ modulo 4. 

\begin{theorem}\label{thm:fermat}
A prime $p\ne2$ is the sum of two squares if and only if $p\equiv1\pmod{4}$.  
\end{theorem}
Theorem \ref{thm:fermat} is usually attributed to Fermat, and appeared in his letter to Mersenne dated Dec~25, 1640 (hence the name \emph{Fermat's Christmas Theorem}), although the statement can already be found in the work of Girard in 1625. The ``only if'' direction is obvious, but the ``if'' direction is far from trivial. Fermat claimed that he had an irrefutable proof, but nobody was able to find the complete proof among his work --- apparently the margin was often too narrow for Fermat. The only clue (in his letters to Pascal and to Digby) is that he used a ``descent argument'': if such a prime $p$ is not of the required form, then one can construct another smaller prime and so on, until a contradiction occurs when one encounters 5, the smallest such prime. More than 100 years later, Euler (1755) gave the first rigorous proof of Theorem \ref{thm:fermat} based on infinite descent. For a detailed history of Theorem \ref{thm:fermat}, see Dickson \cite[Ch. VI, p. 227-231]{Dickson1966}.

\subsection{Which positive integer $n$ can be be written as the sum of two squares?} If $n=x^2+y^2$ ($x,y\in \mathbb{Z}$) and $p\mid n$, then either $p\mid \mathrm{gcd}(x,y)$ or $(x/y)^2\equiv-1\pmod{p}$, and hence either $n/p^2$ is also the sum of two squares or $p\not\equiv3\pmod{4}$ (by the quadratic reciprocity). It follows that each $p\mid n$ with $p\equiv3\pmod{4}$ must appear to an even power. On the other hand, the familiar Diophantus identity $$\left(a^2 + b^2\right)\left(c^2 + d^2\right) =\left(ac-bd\right)^2 + \left(ad+bc\right)^2$$ shows that a product of integers of the form $x^2+y^2$ is also of the same form. Combining with Theorem \ref{thm:fermat} we obtain:

\begin{corollary}
  A positive integer $n$ is of the form $n=x^2+y^2$ if and only if each prime factor $p\equiv3\pmod{4}$ of $n$ appears to an even power.
\end{corollary}

\subsection{In how many different ways can one represent $n$ as the sum of two squares?}

\begin{definition}
To answer this question, we naturally define the \emph{representation number} $$r(n):=\#\{(x,y)\in \mathbb{Z}^2: n=x^2+y^2\}.$$ 
\end{definition}

In particular, $n$ is of the form $x^2+y^2$ if and only if $r(n)>0$.

\begin{example}
  \begin{align*}
    4&= 0^2+(\pm2)^2=(\pm2)^2+0^2, & r(4)=4,\\
  5&=(\pm1)^2+(\pm2)^2=(\pm2)^2+(\pm1)^2, &r(5)=8,\\
 25&=0^2+(\pm5)^2=(\pm3)^2+(\pm4)^2=(\pm4)^2+(\pm3)^2, &r(25)=12. \\
\end{align*}
\end{example}

In his book \emph{Fundamenta nova theoriae functionum ellipticarum} (1829), Jacobi proved the following general formula for the representation numbers.

\begin{theorem}[Jacobi]\label{thm:jacobi}
  $$r(n)=4\left(\sum_{d|n\atop d\equiv1\bmod4}1-\sum_{d|n\atop d\equiv3\bmod4}1\right).$$
\end{theorem}
As a by-product, Jacobi's formula  shows that
  \begin{align*}
    p\equiv1\pmod{4}&\Rightarrow\ r(p)=4(2-0)=8>0, \\
    p\equiv3\pmod{4}&\Rightarrow\ r(p)=4(1-1)=0.
  \end{align*}
which  gives an immediate (and different) proof of Theorem \ref{thm:fermat}!

\subsection{Jacobi's proof} Jacobi's proof of Theorem \ref{thm:jacobi} involves \emph{Jacobi's theta series}, $$\theta:=\sum_{n\in \mathbb{Z}}q^{n^2}=1+2q+2q^4+2q^9+\cdots.$$ The representation numbers $r(n)$ naturally appears as the $n$-th coefficients of the square of Jacobi's theta series $$\theta^2=\left(\sum_{n\in \mathbb{Z}}q^{n^2}\right)^2=\sum_{n\ge0}r(n)q^n=1 + 4 q + 4 q^2 + 4 q^4 + 8 q^5 +4 q^8+\cdots.$$ Jacobi used his theory of elliptic functions (including his famous \emph{Triple Product Identity}) to derive the formula (\cite[p. 107]{Jacobi1829}) $$\theta^2=1+4\left(\frac{q}{1-q}-\frac{q^3}{1-q^3}+\frac{q^5}{1-q^5}-\frac{q^7}{1-q^7}+\cdots\right),$$ which is easily seen to be equivalent to Theorem \ref{thm:jacobi}.

\subsection{Another proof using modular forms} An alternative way of evaluating $\theta^2$ is to view $q=e^{2\pi i\tau}$, and $\theta=\theta(\tau)$ as a holomorphic function on  the upper half plane $$\mathcal{H}:=\{\tau=\sx+i\sy\in \mathbb{C}:\sy=\Im(\tau)>0\}.$$  The function $\theta(\tau)$ satisfies two transformation rules (see \cite[Proposition 9]{Zagier2008})  $$\theta(\tau+1)=\theta(\tau),\quad \theta\left(-\frac{1}{4\tau}\right)=(-2i\tau)^{1/2}\theta(\tau).$$
The first rule is clear by the periodicity of the exponential function. The second rule can be proved using the Poisson summation formula and also plays a key role in Riemann's proof the functional equation of the Riemann zeta function (see \cite[\S 4.9]{Diamond2005}). These rules amounts to say that $$\theta(\tau)\in M_{1/2}(\Gamma_1(4))$$ is a \emph{modular form} of weight $1/2$ and level $\Gamma_1(4)$. Jacobi's theta series  $\theta(\tau)$ and its variants (under the general name of \emph{theta series}) form one of most important classes of modular forms.

It follows that $$\theta^2(\tau)\in M_1(\Gamma_1(4))$$ is a modular form of weight 1 and level $\Gamma_1(4)$. The space $M_1(\Gamma_1(4))$ is in fact 1-dimensional (\cite[Proposition 3]{Zagier2008} or \cite[Theorem 3.6.1]{Diamond2005}), so if one can construct another a modular form of weight 1 and level $\Gamma_1(4)$, then it has to be a scalar multiple of $\theta^2(\tau)$. We next construct such a modular form using \emph{Eisenstein series}, another of most important classes of modular forms.

\begin{definition}
Let $\chi: (\mathbb{Z}/4 \mathbb{Z})^\times\isoarrow \{\pm1\}$ be the (unique) non-trivial character. We define an \emph{Eisenstein series}
\begin{equation}
  \label{eq:eisenstein}
  G_k^{\chi}(\tau)=\sum_{(0,0)\ne (c,d)\in \mathbb{Z}^2, 4|c}\frac{\chi(d)}{(c\tau+d)^k},
\end{equation} where $\chi(d)$ is understood to be 0 when $(4,d)\ne1$.
\end{definition}

 When $k\ge3$, the series (\ref{eq:eisenstein}) is absolutely convergent and is nonzero only when $k$ is odd. When $k\ge3$ is odd, it defines a modular form $G_{k}^\chi(\tau)\in M_k(\Gamma_1(4))$ of weight $k$, level $\Gamma_1(4)$ and character $\chi$.  The constant term of the $q$-expansion of $G_k^\chi(\tau)$ is nonzero and we let $E_k^{\chi}(\tau)$ be a scalar multiple of $G_k^\chi(\tau)$ so the constant term is normalized to be 1.  This \emph{normalized Eisenstein series} $E_k^\chi(\tau)$ then have the explicit $q$-expansion  (see \cite[\S 4.5]{Diamond2005})
 \begin{equation}
   \label{eq:qexpansion}
   E_k^{\chi}(\tau)=1+c_k^\chi\cdot\sum_{n\ge1}\left(\sum_{d|n}\chi(d)d^{k-1}\right)q^n,
 \end{equation}
 where $c_k^\chi=2/L(1-k,\chi)$ is related to a special value of the Dirichlet $L$-function $L(s,\chi)$.

When $k=1$, the series (\ref{eq:eisenstein}) is not absolutely convergent, but one can still  suitably modify it  to obtain a modular form $$E_1^\chi(\tau)\in M_1(\Gamma_1(4))$$ with the same formula (\ref{eq:qexpansion}) for its $q$-expansion, either using  the Weierstrass $\sigma$-function (see \cite[\S 4.8]{Diamond2005}), or using the analytic continuation of
\begin{equation}
  \label{eq:G1}
  G_1^{\chi}(\tau,s):=\sum_{(0,0)\ne (c,d)\in \mathbb{Z}^2, 4|c}\frac{\chi(d)}{(c\tau+d)|c\tau+d|^{2s}},\quad \Re(s)>1/2 
\end{equation}
to $s=0$ (see \cite[\S 7.2]{Miy89}). In particular, when $k=1$ the formula (\ref{eq:qexpansion}) simplifies to $$E_1^\chi(\tau)=1+c_1^\chi\cdot\sum_{n\ge1}\left(\sum_{d|n}\chi(d)\right)q^n.$$

 As explained, $E_1^\chi(\tau)$ must be a scalar multiple of $\theta^2(\tau)$. Since both of them have constant coefficient 1, we indeed have the equality \begin{equation}
   \label{eq:Jacobi}
   \theta^2(\tau)=E_1^\chi(\tau).
 \end{equation}  Comparing the coefficient before $q$, we obtain $c_1^\chi=4$ and hence $$r(n)=4\sum_{d|n}\chi(d),$$ which proves Theorem \ref{thm:jacobi}.

\begin{remark}
    As a by-product of the proof, we also obtain $L(0,\chi)=\frac{1}{2}$ from $c_1^\chi=4$. Via the functional equation of $L(s,\chi)$, this is equivalent to the famous Leibniz formula for $\pi$ (1676), $$L(1,\chi)=1-\frac{1}{3}+\frac{1}{5}-\frac{1}{7}+\cdots=\frac{\pi}{4}.$$
\end{remark}

 To summarize, Jacobi's Theorem \ref{thm:jacobi} can be proved using the identity of two modular forms (\ref{eq:Jacobi}), namely using a relation of the form 
 \begin{center}
\fbox{theta series $\longleftrightarrow$ Eisenstein series.}
 \end{center} Notice that the Fourier coefficients of theta series encodes representation numbers of quadratic forms, while the Fourier coefficients of Eisenstein series are generalized divisor sums which are more explicit.

 \section{Siegel--Weil formula}

\subsection{Siegel's formula} 
Siegel \cite{Sie35} generalizes the formula (\ref{eq:Jacobi}) from the binary quadratic form $x^2+y^2$ to more general quadratic forms in arbitrary number of variables. Let $\Lambda$ be a positive definite quadratic lattice over $\mathbb{Z}$ of rank $m$ with quadratic form $Q: \Lambda \rightarrow \mathbb{Z}$. Denote by $(\ ,\ )$ the associated symmetric bilinear form, defined by $$(x,y):=Q(x+y)-Q(x)-Q(y),$$ (so $Q(x)=\frac{1}{2}(x,x)$). Denote by $\Sym_n(\mathbb{Z})$ the set of symmetric matrices whose diagonal entries are in $\mathbb{Z}$ and whose off-diagonal entries are in $\frac{1}{2}\mathbb{Z}$. Denote by $\Sym_n(\mathbb{Z})_{\ge0}\subseteq\Sym_n(\mathbb{Z})$ the subset of positive semi-definite matrices.

\begin{definition}
For $T\in\Sym_n(\mathbb{Z})_{\ge0}$, define the (generalized) \emph{representation number} $$r_\Lambda(T):=\#\{(x_1,\ldots,x_n)\in \Lambda^n: {\textstyle\frac{1}{2}}((x_i,x_j))_{1\le i,j\le n}=T\}.$$ Define \emph{Siegel's theta series}
\begin{equation}
  \label{eq:siegeltheta}
  \theta_\Lambda(\sz):=\sum_{T\in \Sym_n(\mathbb{Z})_{\ge0}}r_\Lambda(T)q^T, \quad q^T:=e^{2\pi i\tr T\sz},
\end{equation}
 a holomorphic function on \emph{Siegel's half space} $$\mathcal{H}_n:=\{\sz=\sx+i \sy: \sx\in \Sym_n(\mathbb{R}), \sy\in \Sym_n(\mathbb{R})_{>0}\}.$$ Using the Poisson summation formula, Siegel proved that $\theta_\Lambda(\sz)$ is a Siegel modular form on $\mathcal{H}_n$ of weight $m/2$.
\end{definition}

\begin{example}
Notice that when $n=1$, Siegel's half space $\mathcal{H}_n$ recovers the usual upper half plane $\mathcal{H}$, and Siegel's theta series $\theta_\Lambda(\sz)$ recovers Jacobi's theta series  $$\theta_\Lambda(\tau):=\sum_{n\ge0}r_\Lambda(n)q^n,\quad r_\Lambda(n):=\{x\in \Lambda: Q(x)=n\}.$$
\end{example}

In general, the theta series $\theta_\Lambda$ for a lattice $\Lambda$ may fail to be an Eisenstein series on the nose, but Siegel's formula shows that the weighted average of theta series within its genus class is always a Siegel Eisenstein series on $\mathcal{H}_n$:
 \begin{center}
\fbox{weighted average of theta series  $\longleftrightarrow$ Siegel Eisenstein series.}
\end{center}

More precisely, recall that  two quadratic lattices $\Lambda,\Lambda'$ are in the same \emph{genus}, if they are isomorphic over  $\mathbb{R}$ and over $\mathbb{Z}_p$ for all primes $p$. Denote by $\Gen(\Lambda)$ the set of isomorphism classes of quadratic lattices in the same genus of $\Lambda$. Denote by $\Aut(\Lambda)$ the automorphism group of $\Lambda$ as a quadratic lattice.

\begin{theorem}[Siegel] The following identity holds:
\begin{equation}
  \label{eq:Siegel}
\frac{
  \sum\limits_{\Lambda'\in\Gen(\Lambda)}\frac{1}{\#\Aut(\Lambda')}\cdot\theta_{\Lambda'}(\sz)}{ \sum\limits_{\Lambda'\in\Gen(\Lambda)}\frac{1}{\#\Aut(\Lambda')}}
=E_\Lambda(\sz).
\end{equation}
 Here $E_\Lambda(\sz)$ is a certain normalized Siegel Eisenstein series on $\mathcal{H}_n$ of weight $m/2$.
\end{theorem}

\begin{example}\label{exa:jacobisiegel}
 Consider the case $m=2$, $n=1$ and $\Lambda=\mathbb{Z}^2$ equipped with the quadratic form $Q=x^2+y^2$. Then $$\theta_\Lambda(\tau)=\theta^2(\tau),\quad E_\Lambda(\tau)= E_1^\chi(\tau).$$ In this case  $\Gen(\Lambda)$ is a singleton and Siegel's formula recovers (\ref{eq:Jacobi}).
\end{example}

\begin{example}[{cf. \cite[V.2.3]{Ser73}}]
  Siegel's formula is extremely useful in studying the arithmetic of quadratic forms. For example, one can deduce his famous \emph{mass formula} (also known as the Smith--Minkowski--Siegel mass formula), which computes the \emph{mass} of $\Gen(\Lambda)$, defined to be weighted size $$\sum_{\Lambda'\in\Gen(\Lambda)}\frac{1}{\#\Aut(\Lambda')}$$ as an Euler product of local factors indexed by primes $p$.

  For example, consider the simplest case when $\Lambda$ is
  \begin{itemize}
  \item \emph{unimodular}, i.e., $\det(\frac{1}{2}((x_i,x_j)_{i,j=1}^n)\in\{\pm1\}$ for a $\mathbb{Z}$-basis $\{x_1,\ldots,x_n\}$ of $\Lambda$,
  \item and \emph{even}, i.e., 2 divides $Q(x)$ for all $x\in \Lambda$.
  \end{itemize}
The rank $m$ of any unimodular even lattice $\Lambda$ is necessary a multiple of 8. Siegel's mass formula computes the mass of $\Gen(\Lambda)$ explicitly as
  \begin{align*}
    \sum_{\Lambda'\in\Gen(\Lambda)}\frac{1}{\#\Aut(\Lambda')}&=\frac{2\zeta (2)\zeta (4)\cdots \zeta (m-2)\zeta(m/2)}{\vol(S^0)\vol(S^1)\cdots \vol(S^{m-1})}
    =
    {B_{{m/2}} \over m}\prod _{{1\leq j<m/2}}{B_{{2j}} \over 4j}.
  \end{align*}
where $\zeta(s)$ is the Riemann zeta function, $\vol(S^{k-1})=\frac{2\pi^{k/2}}{\Gamma(k/2)}$ is the volume of the unit $(k-1)$-sphere and $B_k$ is the $k$-th Bernoulli number.
\end{example}
\begin{example}[{cf. \cite[VII.6.6]{Ser73}}]
  Let $\Lambda$ be the root lattice of type $E_8$, defined by $$\Lambda:=\{x=(x_1,\ldots, x_8)\in \mathbb{Z}^8\cup (\mathbb{Z}+1/2)^8: \sum_{i=1}^8x_i\equiv0\pmod{2}\},\quad Q=\sum_{i=1}^8 x_i^2.$$ Then $\Lambda$ is unimodular and even. Siegel's mass formula computes that $$\sum_{\Lambda'\in\Gen(\Lambda)}\frac{1}{\#\Aut(\Lambda')}={B_{4} \over 8}{B_{2} \over 4}{B_{4} \over 8}{B_{6} \over 12}={-1/30 \over 8}\;{1/6 \over 4}\;{-1/30 \over 8}\;{1/42 \over 12}={1 \over 696729600}.$$ The fact that $\#\Aut(\Lambda)=696729600$ (which is also the order of the Weyl group of type $E_8$) then implies that $\Lambda$ is the \emph{unique} unimodular and even lattice of rank 8.

In this case $E_\Lambda(\tau)$ is related to the classical Eisenstein series $E_4(\tau)\in M_4(\SL_2(\mathbb{Z}))$ of weight 4 and level 1 by
\begin{align*} 
  E_\Lambda(\tau)=E_4(2\tau)&=1+240\sum_{n\ge1}\left(\sum_{d|n}d^3\right)q^{2n}
  =1+240q^2+2160q^4+6720q^6+\cdots. 
\end{align*}
 Here the factor $240$ agrees with the constant $c_k=2/\zeta(1-k)$ for the weight $k=4$.  Siegel's formula then implies that $$r_\Lambda(2n):=\#\left\{x=(x_i)\in \Lambda: \sum_{i=1}^8 x_i^2=2n\right\}=240\sum_{d|n}d^3.$$ In particular, we recover that there are $r_\Lambda(2)=240$ roots in the $E_8$ root system (see Figure \ref{fig:E8}).
\begin{figure}[h]
  \centering
  \includegraphics[scale=0.1]{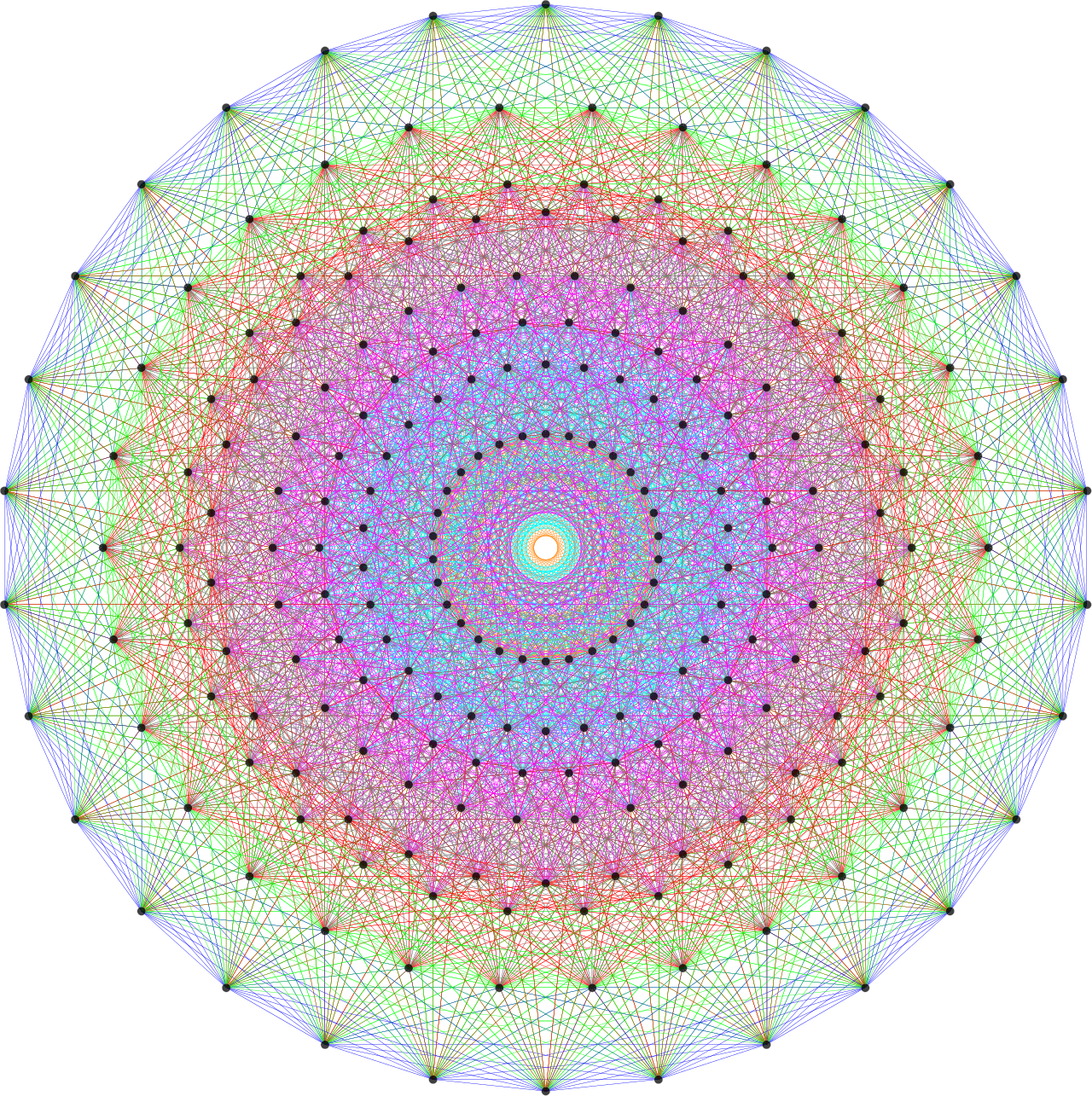} 
  \caption{A 2-dimensional projection of the convex hull of the  $E_8$ root system. The black dots are the 240 roots. The plane the roots are projected into is one of the four planes in $\mathbb{R}^8$ stabilized by a Coxeter element $w$ of the Weyl group. The Coxeter element $w$ acts on this plane via rotation through an angle of $\frac{2\pi}{h}$, where $h=30$ is the Coxeter number. It acts on the 240 roots regularly and divides them into 8 orbits of size 30 around circles of different radii. Colors of 6720 edge lines are varied based on their projected lengths to increase contrast.\protect\footnotemark}
  \label{fig:E8}
\end{figure}
\end{example}

\subsection{Siegel--Weil formula} In a series of works (\cite{Sie36, Sie37, Siegel1951, Sie52}), Siegel further generalized his formula form definite to indefinite quadratic forms and from the base field $\mathbb{Q}$ to totally real fields. In the indefinite case the formula is more difficult to state, as the theta series $\theta_\Lambda(\sz)$ is divergent and it is necessarily to introduce extra weight functions in the definition to ensure convergence. The situation was greatly clarified by Weil \cite{Weil1965} using powerful tools from representation theory, especially due to his use of the \emph{Weil representation}.

\footnotetext{Computer generated picture by Jgmoxness, downloaded from \url{https://commons.wikimedia.org/wiki/File:E8Petrie.svg}}

Let $V$ be a quadratic space over $\mathbb{Q}$ of dimension $m$ with bilinear form $(\ ,\ )$. For simplicity we assume that $m$ is \emph{even} (so the weight $m/2$ of the relevant Siegel modular forms is integral, see Remark \ref{rem:oddweil}).  Consider the reductive dual pair $(G, H)=(\Sp(2n), \OO(V))$, where $\Sp(2n)$ is the symplectic group of the standard $2n$-dimensional symplectic space over $\mathbb{Q}$ and $\O(V)$ is the orthogonal group of $V$.  Let $P=MN\subseteq G=\Sp(2n)$ be the standard Siegel parabolic subgroup, so that under the standard basis we have 
\begin{align*}
M(\mathbb{Q})&=\left\{m(a)=\begin{pmatrix}a & 0\\0 &{}^ta^{-1}\end{pmatrix}: a\in \GL_n(\mathbb{Q})\right\},\\
  N(\mathbb{Q}) &= \left\{n(b)=\begin{pmatrix} 1_n & b \\0 & 1_n\end{pmatrix}: b\in \Sym_n(\mathbb{Q})\right\}.  
\end{align*}

\begin{definition}
Let $\mathbb{A}$  be the ring of ad\`eles of $\mathbb{Q}$. We fix the standard additive character $\psi: \mathbb{A}\rightarrow \mathbb{C}^\times$ whose archimedean component is given by  $\psi_\infty: \mathbb{A}_\infty=\mathbb{R}\rightarrow \mathbb{C}^\times, x\mapsto e^{2\pi i x}$. The \emph{(Schr\"odinger model of the) Weil representation} $\omega=\omega_{V,\psi}$ is the representation of $G(\mathbb{A})\times H(\mathbb{A})$ on the space of Schwartz functions $\sS(V(\mathbb{A})^n)$ such that for any $\varphi\in \sS(V(\mathbb{A})^n)$ and $\mathbf{x}\in V(\mathbb{A})^n$,
\begin{align}\label{eq:weil}
  \begin{aligned}
      \omega(m(a))\varphi(\mathbf{x})&=\chi_V(\det a)|\det a|^{m/2}\varphi(\mathbf{x}\cdot a),\quad &m(a)\in M(\mathbb{A}),\\
\omega(n(b))\varphi(\mathbf{x})&=\textstyle\psi(\frac{1}{2}\tr b\,(\mathbf{x},\mathbf{x}))\varphi(\mathbf{x}),&n(b)\in N(\mathbb{A}),\\
\omega(w)\varphi(\mathbf{x})&=
\widehat \varphi(\mathbf{x}),&w=\left(\begin{smallmatrix}
0  & 1_n\\
  -1_n & 0\\
\end{smallmatrix}\right),\\
\omega(h)\varphi(\mathbf{x})&=\varphi(h^{-1}\cdot\mathbf{x}),& h\in H(\mathbb{A}).
  \end{aligned}
\end{align} 
Here
\begin{itemize}
\item $\chi_V: \mathbb{A}^\times/\mathbb{Q}^\times\rightarrow \mathbb{C}^\times$ is the quadratic character corresponds to the quadratic extension $\mathbb{Q}(\sqrt{\disc(V)})/\mathbb{Q}$, and $\disc(V)\in \mathbb{Q}^\times/(\mathbb{Q}^\times)^2$ is the discriminant of $V$ defined to be $$\disc(V):=(-1)^{m \choose 2}\det ((x_i,x_j))_{i,j=1}^m$$ for any $\mathbb{Q}$-basis $\{x_1,\ldots,x_m\}$ of $V$.
\item $|\cdot |: \mathbb{A}^\times\rightarrow \mathbb{R}_{>0}$ is the normalized absolute value.
\item  $(\mathbf{x},\mathbf{y}):=((x_i,y_j))_{i,j=1}^n\in \mathrm{Mat}_n(\mathbb{A})$ for $\mathbf{x}=(x_1,\ldots,x_n)\in V(\mathbb{A})^n$ and $\mathbf{y}=(y_1,\ldots, y_n)\in V(\mathbb{A})^n$,
\item $\widehat\varphi$ is the Fourier transform of $\varphi$ using the self-dual Haar measure on $V(\mathbb{A})^n$ with respect to $\psi$, $$\widehat\varphi(\mathbf{x}):=\int_{V(\mathbb{A})^n}\varphi(\mathbf{y})\psi(\tr(\mathbf{x},\mathbf{y}))\rd \mathbf{y}.$$
\end{itemize} 
\end{definition}

\begin{remark}
  There are many nice expository articles on the Weil representation, see e.g. \cite{Kud96, Pra93, Pra98}. See also a more general formula of the Weil representation in \cite[Proposition 4.3]{Kud96}). 
\end{remark}

\begin{example}
  When $n=1$ we have $G=\Sp(2)=\SL(2)$. The standard Siegel parabolic is the standard Borel subgroup of $\SL(2)$ consisting of upper triangle matrices $P=\left\{\left(\begin{smallmatrix}a & b \\ 0 & a^{-1}\end{smallmatrix}\right)\right\}$, and $M=\left\{\left(\begin{smallmatrix}a & 0 \\ 0 & a^{-1}\end{smallmatrix}\right)\right\}$, $N=\left\{\left(\begin{smallmatrix}1 & b \\ 0 & 1\end{smallmatrix}\right)\right\}$ are the diagonal and upper unipotent matrices respectively. In this case the first two formulas in (\ref{eq:weil}) simplify to
  \begin{equation}
    \label{eq:weilSL2}
    \omega\left(\begin{smallmatrix}a & 0 \\ 0 & a^{-1}\end{smallmatrix}\right)\varphi(x)=\chi_V(a)|a|^{m/2}\varphi(xa),\quad \omega\left(\begin{smallmatrix}1 & b \\ 0 & 1\end{smallmatrix}\right)\varphi(x)=\textstyle\psi(\frac{1}{2}b(x,x))\varphi(x),
  \end{equation}
 for any $a\in \mathbb{A}^\times$, $b\in \mathbb{A}$, $\varphi\in \sS(V(\mathbb{A}))$ and $x\in V(\mathbb{A})$. 
\end{example}

Our next goal is to use the Weil representation $\omega$ to construct theta series and Siegel Eisenstein series, starting from any common choice of a Schwartz function $\varphi\in \sS(V(\mathbb{A})^n)$.

\begin{definition}
  Associated to $\varphi\in \sS(V(\mathbb{A})^n)$, define the (two-variable) \emph{theta function} $$\theta(g,h,\varphi):=\sum_{\mathbf{x}\in V^n}\omega(g)\varphi(h^{-1}\mathbf{x}),\quad g\in G(\mathbb{A}), h\in H(\mathbb{A}).$$ Then $\theta(g,h,\varphi)$ is invariant under $G(\mathbb{Q})\times H(\mathbb{Q})$ (automorphic on both $G$ and $H$ in a broad sense).
\end{definition}

\begin{remark}\label{rem:thetalifting}
 Using $\theta(g,h,\varphi)$ as an integral kernel allows one to lift automorphic forms on $G$ to automorphic forms on $H$ (and vice versa):  for a cuspidal automorphic representation $\pi$ of $G(\mathbb{A})$ and $\phi\in \pi$, define the \emph{theta lift} $\theta_\varphi(\phi)$ of $\phi$ to $H(\mathbb{A})$ by the Petersson inner product on $G(\mathbb{Q})\backslash G(\mathbb{A})$, $$\theta_\varphi(\phi)(h):=\langle \theta(-,h,\varphi),\phi\rangle_\mathrm{Pet}=\int_{G(\mathbb{Q})\backslash G(\mathbb{A})}\theta(g,h,\varphi)\overline{\phi(g)}\rd g.$$  Then $\theta_\varphi(\phi)$ is an automorphic form on $H(\mathbb{A})$. This may be viewed as the starting point of the modern theory of \emph{theta correspondence}, which is indispensable in the study of automorphic forms and the Langlands correspondence. We refer to Gan \cite{Gan14} for an excellent recent survey on theta correspondence.
\end{remark}

  \begin{example}\label{exa:convertclassical}
Assume that $V$ is positive definite, then the \emph{theta integral}
\begin{equation}
  \label{eq:thetaintegral}
  \int_{H(\mathbb{Q})\backslash H(\mathbb{A})}\theta(g,h,\varphi)\ \rd h
\end{equation}
 or in other words, the theta lift of the constant function $\mathbf{1}$ on $H(\mathbb{A})$ to $G(\mathbb{A})$, is closely related to Siegel's theta series. More precisely, for $\Lambda\subseteq V$ a lattice over $\mathbb{Z}$, we take the Schwartz function $\varphi=(\otimes_p\varphi_p)\otimes \varphi_\infty\in \sS(V(\mathbb{A})^n)$ such that
\begin{itemize}
\item $\varphi_p\in\sS(V(\mathbb{Q}_p)^n)$ is the characteristic function of $(\Lambda \otimes \mathbb{Z}_p)^n$,
\item $\varphi_\infty\in \sS(V(\mathbb{R})^n)$ is the standard Gaussian function $\varphi_\infty(\mathbf{x})=e^{-\pi \tr (\mathbf{x},\mathbf{x})}$.
\end{itemize}

For $\sz=\sx+i\sy\in \mathcal{H}_n$, we consider $g_\sz=n(\sx)m(a)\in G(\mathbb{R})$, where $a\in\GL_n(\mathbb{R})$ such that $\sy=a{}^{t} a$. Then $g_\sz\cdot i 1_n=\sz$.
By (\ref{eq:weil}) we have
\begin{equation}
  \label{eq:archimedean}
  \omega_\infty(g_\sz)\varphi_\infty(\mathbf{x})=\chi_\infty(\det a) |\det a|^{m/2}\cdot q^{\frac{1}{2}(\mathbf{x},\mathbf{x})}.
\end{equation}
 Define the \emph{classical theta integral}
\begin{equation}
  \label{eq:classicalthetaintegral}
  \chi_\infty(\det a)^{-1}|\det a|^{-m/2}\cdot \int_{H(\mathbb{Q})\backslash H(\mathbb{A})}\theta(g_{\sz},h,\varphi)\ \rd h.
\end{equation}
Then it recovers the weighted average of theta series in (\ref{eq:Siegel}) (see e.g., \cite[\S4.6]{Han13}, \cite[\S7]{Kudla2014}). 

In fact, let $\KL\subseteq H(\mathbb{A})$ be the stabilizer of $\Lambda$, then we have a bijection $$H(\mathbb{Q})\backslash H(\mathbb{A})/\KL\isoarrow \Gen(\Lambda),\quad h \mapsto \Lambda':=h(\Lambda \otimes \widehat{\mathbb{Z}})\cap V.$$ Let $\{h_i\}$ be a complete set of representatives of $H(\mathbb{Q})\backslash H(\mathbb{A})/\KL$ and let $\{\Lambda_i\}$ be the corresponding representatives of $\Gen(\Lambda)$ under this bijection.  Then $$\int_{H(\mathbb{Q})\backslash H(\mathbb{A})}\theta(g_{\sz},h,\varphi)\ \rd h=\sum_i\int_{H(\mathbb{Q})\backslash H(\mathbb{Q})h_i \KL}\theta(g_{\sz},h,\varphi)\rd h.$$ Using $H(\mathbb{Q})\cap h_i \KL h_i^{-1}=\Aut(\Lambda_i)$, each summand evaluates to $$\int_{H(\mathbb{Q})\backslash H(\mathbb{Q})h_i \KL h_i^{-1}} \theta(g_{\sz}, h h_i,\varphi)\rd h=\frac{1}{\#\Aut(\Lambda_i)}\int_{h_i\KL h_i^{-1}}\theta(g_{\sz},h h_i, \varphi)\rd h.$$ Unfolding the definition, the second integral equals $$\int_{\KL}\theta(g_{\sz},  h_i h, \varphi)\rd h=\int_K\sum_{\mathbf{x}\in V^n}\omega_\infty(g_{\sz})\varphi(h^{-1}h_i^{-1}\mathbf{x})\rd h,$$  which by our choice of $\{\varphi_p\}$ evaluates to $$\vol(K)\sum_{\mathbf{x}\in \Lambda_i^n}\omega_\infty(g_{\sz})\varphi_\infty(\mathbf{x}).$$ Thus combining with (\ref{eq:archimedean}) we know that the classical theta integral (\ref{eq:classicalthetaintegral}) evaluates to $$\vol(K)\sum_i \frac{1}{\#\Aut(\Lambda_i)}\cdot\sum_{\mathbf{x}\in \Lambda_i^n}q^{\frac{1}{2}(\mathbf{x},\mathbf{x})}=\vol(K)\sum_i \frac{1}{\#\Aut(\Lambda_i)}\cdot\theta_{\Lambda_i}(\sz).$$ Finally notice that $$\vol(K)\sum_i \frac{1}{\#\Aut(\Lambda_i)}=\vol(H(\mathbb{Q})\backslash H(\mathbb{A})),$$ thus if we normalize the Haar measure $\rd h$ such that $\vol(H(\mathbb{Q})\backslash H(\mathbb{A}))=1$, then  the classical theta integral (\ref{eq:classicalthetaintegral}) recovers the weighted average of theta series in (\ref{eq:Siegel}).

\end{example}

\begin{definition}
  Associated to $\varphi\in \sS(V(\mathbb{A})^n)$, also define the \emph{Siegel Eisenstein  series} $$E(g, s,\varphi):=\sum_{\gamma\in P(\mathbb{Q})\backslash G(\mathbb{Q})}\Phi_\varphi(\gamma g,s),\quad g\in G(\mathbb{A}), s\in \mathbb{C},$$ where
  \begin{equation}
    \label{eq:SWsection}
    \sS(V(\mathbb{A})^n)\rightarrow \Ind_{P(\mathbb{A})}^{G(\mathbb{A})}(\chi_V|\cdot|^{s+\frac{n+1}{2}}),\quad \varphi\mapsto \Phi_\varphi(g,s):=\omega(g)\varphi(0)\cdot |\det a(g)|^{s-s_0}
  \end{equation}
 is the \emph{standard Siegel--Weil section} of the degenerate principal series representation of $G(\mathbb{A})$ and $$s_0:=\frac{m-(n+1)}{2}.$$ Here we write $g=nm(a)k$ under the Iwasawa decomposition $G(\mathbb{A})=N(\mathbb{A})M(\mathbb{A})K$ for $K$ the standard maximal open compact subgroup of $G(\mathbb{A})$, and the quantity $|\det a(g)|:=|\det a|$ is well-defined.
\end{definition}

\begin{example}\label{exa:classicaleisenstein}
  Similarly define the \emph{classical Siegel Eisenstein series} $$E(\sz,s,\varphi)\coloneqq \chi_\infty(\det a)^{-1}|\det a|^{-m/2}\cdot E(g_{\sz},s, \varphi).$$  When $V$ is positive definite and $\varphi$ is chosen as in Example \ref{exa:convertclassical}, the special value $E(\sz, s_0,\varphi)$ at $s=s_0$ essentially recovers $E_\Lambda(\sz)$ in (\ref{eq:Siegel}).

  For example, consider the case $m=2$, $n=1$ and $\Lambda=\mathbb{Z}^2$ equipped with the quadratic form $Q=x^2+y^2$. We have $\disc(V)=-1$ and the quadratic character $\chi_V$ corresponds to the quadratic extension  $\mathbb{Q}(i)/\mathbb{Q}$, and hence corresponds the Dirichlet character $\chi: (\mathbb{Z} /4 \mathbb{Z})^\times\simeq\{\pm1\}$.  For $\tau=\sx+i\sy\in \mathcal{H}$, we have $$g_\tau=\left(\begin{smallmatrix} 1 & \sx \\ 0 & 1\end{smallmatrix}\right)\left(\begin{smallmatrix}\sy^{1/2} & 0 \\ 0 & \sy^{-1/2}\end{smallmatrix}\right)\in\SL_2(\mathbb{R}).$$  For $\gamma= \left(\begin{smallmatrix}a & b\\ c & d\end{smallmatrix}\right)\in \SL_2(\mathbb{Q})$, one can compute that the Siegel--Weil section evaluates to $$\Phi_\varphi(\gamma g_\tau,s)=\frac{\chi(d)\cdot\sy^{1/2}}{c\tau+d}\cdot \Im(\gamma \tau)^{s/2}=\frac{\chi(d)}{c\tau+d}\cdot \frac{\sy^{(1+s)/2}}{|c\tau+d|^s},$$ using (\ref{eq:archimedean}) together with (\ref{eq:weilSL2}) (or the more general \cite[Proposition 4.3]{Kud96}).  Comparing with (\ref{eq:G1}) we see that the classical Eisenstein series $$E(\tau, 0, \varphi)=\sy^{-1/2}\cdot E(g_\tau, 0,\varphi)$$ at $s_0=0$ essentially recovers  the Eisenstein series $E_1^\chi(\tau)$ of weight 1 (up to a nonzero constant).
\end{example}

 The Siegel Eisenstein series $E(g,s,\varphi)$ converges absolutely when $\Re(s)>\frac{n+1}{2}$. It has a meromorphic continuation to $s\in \mathbb{C}$ and satisfies a functional equation relating $s\leftrightarrow -s$ (i.e., centered at $s=0$). The Siegel--Weil formula gives a precise identity of the form 
 \begin{center}
\fbox{theta integral $\longleftrightarrow$ value of Siegel Eisenstein series at $s=s_0$.}
\end{center}
 Notice that $s=s_0$ is the unique point such that the map $\varphi\mapsto \Phi_\varphi(g,s)$ in (\ref{eq:SWsection}) is $G(\mathbb{A})$-equivariant, so that both sides of identity at least have the same transformation behavior with respect to the Weil representation.

\begin{theorem}[Siegel--Weil formula \cite{Weil1965,KR88, KR88a}]\label{thm:SWF} Let $\alpha$ be the dimension of a maximal isotropic subspace of $V$. If $\alpha=0$ (i.e., $V$ is anisotropic) or $\alpha>0$ and $m-\alpha>n+1$, then $E(g,s,\varphi)$ is holomorphic at $s_0$ and $$\kappa\cdot\int_{H(\mathbb{Q})\backslash H(\mathbb{A})}\theta(g,h,\varphi)\ \rd h= E(g,s_0,\varphi)
        ,$$ where $\kappa=1$ if $m>n+1$ or $\kappa=2$ otherwise. Here the Haar measure $\rd h$ is normalized so that $\vol(H(\mathbb{Q})\backslash H(\mathbb{A}))=1$.
      \end{theorem}

\begin{example}
When $V$ is positive definite (so $\alpha=0$), the Siegel--Weil formula recovers Siegel's formula (see Examples \ref{exa:classicaleisenstein} and \ref{exa:convertclassical}), and in particular recovers Jacobi's formula  as a special case of $m=2$, $n=1$, $(G,H)=(\Sp(2),\O(2,0))$
(see Example \ref{exa:jacobisiegel}).
\end{example}

\begin{remark}
The condition in Theorem \ref{thm:SWF} is known as \emph{Weil's convergence condition}, which ensures the the convergence of the theta integral. It is a long effort starting with the work of Kudla--Rallis \cite{KR94} to generalize the Siegel--Weil formula outside the convergence range and for all reductive dual pairs of classical groups. We refer to Gan--Qiu--Takeda \cite{GQT14} for the most general Siegel--Weil formula and a nice summary of the literature and history.
\end{remark}

\begin{remark}\label{rem:oddweil}
  When $m$ is odd, the theta series and Eisenstein series have half integral weights $m/2$, and are automorphic forms on the metaplectic cover $\tilde G=\Mp(2n)$ of $\Sp(2n)$. In this case the Weil representation needs to be modified to be a representation of $\tilde G(\mathbb{A})\times H(\mathbb{A})$ and the Siegel--Weil formula still holds after modification.
\end{remark}

\section{Geometric Siegel--Weil formula}

In this section we discuss an example of the Siegel--Weil formula in the indefinite case $(G,H)=(\Sp(4),\O(2,2))$ originating from the classical work of Hurwitz, and use it to motivate the more general geometric Siegel--Weil formula.

\subsection{Hurwitz class number relation} \label{sec:hurwitz-class-number}

\begin{definition}
For any positive integer $D$, the \emph{Hurwitz class number} $H(D)$ is defined to be the weighted size of $\SL_2(\mathbb{Z})$-equivalence classes of positive
  definite binary quadratic forms
\begin{center}
    $ax^2+bxy+cy^2$ with
  discriminant $b^2-4ac=-D$,\quad{ $a,b,c\in \mathbb{Z}$}.
\end{center}
Here the forms equivalent to $a(x^2+y^2)$ ($D=4a^2$) and $a(x^2+xy+y^2)$ ($D=3a^2$) are counted with multiplicities $1/2$ and $1/3$ respectively, due to extra symmetry.
\end{definition}

\begin{example}
  \begin{align*}
    H(3)  =\frac{1}{3}& \longleftrightarrow\{x^2+xy+y^2\}\\
    H(8) = 1 & \longleftrightarrow\{x^2+2y^2\}\\
    H(11) = 1& \longleftrightarrow\{x^2+xy+3y^2\}\\
 H(12)  =\frac{4}{3}& \longleftrightarrow \{2(x^2+xy+y^2), x^2+3y^2\}\\
  H(20)= 2 &\longleftrightarrow \{x^2+5y^2, 2x^2+2xy+3y^2\}.\\
\end{align*}
\end{example}

\begin{example}
   Table \ref{tab:hurwitz} lists the first few Hurwitz class numbers.
  \begin{table*}[h]
  \centering
  \begin{tabular}[h]{|c|c|c|c|c|c|c|c|c|c|c|c|c|}
    $D$ & 3 & 4 & 7 & 8 & 11 & 12 & 15 & 16 & 19 & 20 & 23 & 24\\\hline
    $H(D)$ & $\frac{1}{3}$ & $\frac{1}{2}$ & 1 & 1 & 1 & $\frac{4}{3}$ & 2 & $\frac{3}{2}$ & 1 & 2 &3 & 2 \\
  \end{tabular}
  \caption{Hurwitz class numbers}
  \label{tab:hurwitz}
\end{table*}
\end{example}

  \begin{example}
When $-D$ is a fundamental discriminant and $D>4$, the Hurwitz class number $H(D)$ is equal to the class number of the imaginary quadratic field $\mathbb{Q}(\sqrt{-D})$.
  \end{example}

Understanding these class numbers $H(D)$ remains a central subject in algebraic number theory. The following remarkable formula, which we call the \emph{Hurwitz class number relation} or the \emph{Hurwitz formula}, gives an elementary expression for a certain sum of Hurwitz class numbers.

\begin{theorem}[Kronecker \cite{Kro60}, Gierster \cite{Gie83}, Hurwitz \cite{Hur85}]\label{thm:hurwitz}
If $m$ is not a perfect square, then
\begin{equation}\label{eq:hurwitz}
\sum_{d d'=m} \max\{ d, d'\}=\sum_{t\in \mathbb{Z}, 4m-t^2>0}H(4m-t^2).  
\end{equation}  
\end{theorem}

\begin{example}
When $m=3$, the Hurwitz class number relation says
\begin{align*}
  3+3 & = H(3)+H(8)+H(11)+H(12)+H(11)+H(8)+H(3)  \\
  & =\frac{1}{3}+1+1+\frac{4}{3}+1+1+\frac{1}{3}.
\end{align*} 
When $m=5$, the Hurwitz class number relation says 
\begin{align*}
  5+5 & = H(4)+H(11)+H(16)+H(19)+H(20)+H(19)+H(16)+H(11)+H(4)  \\
   & =\frac{1}{2}+1+\frac{3}{2}+1+2+1+\frac{3}{2}+1+\frac{1}{2},
 \end{align*}
 A quite nontrivial way to decompose the integer 6 and 10 respectively!
\end{example}

Hurwitz \cite{Hur85} proved this formula using the modular $j$-invariant and the modular polynomial $\Phi_m(x, y)\in \mathbb{Z}[x, y]$ of level $m$ (which defines the modular curve $Y_0(m)$ over $\mathbb{Z}$). He recognized the LHS of (\ref{eq:hurwitz}) as $\deg \Phi_m(x,x)$ and computed this degree in a different way involving Hurwitz class numbers to arrive at the formula (\ref{eq:hurwitz}).

\subsection{A geometric proof}(cf. {\cite{Gross1993}}\label{sec:geometric-proof}). 
From the modern point of view, we have a nice geometric interpretation  of Hurwitz's proof, in terms of the geometry of the modular curve $$Y(\mathbb{C})=\SL_2(\mathbb{Z})\backslash \mathcal{H}.$$ The modular curve $Y(\mathbb{C})$ is the moduli space of elliptic curves:
\begin{align*} Y&=\{E: \text{elliptic curve over }\mathbb{C} \text{ (up to isomorphism)}\}\\
 \tau&\mapsto E_\tau=\mathbb{C}/(\mathbb{Z}+\mathbb{Z}\tau).
    \end{align*}  which allows one to define a canonical model of $Y$ as an algebraic curve over $\mathbb{Q}$.     Each elliptic curve $E_\tau$ has a Weierstrass equation  $$E_\tau: y^2=x^3+A_\tau x+B_\tau.$$  The $j$-invariant $$j(E_\tau):=1728\cdot \frac{4A_\tau^3}{4A_\tau^3+27B_\tau^2}$$ only depends on the isomorphism class of $E_\tau$ and gives rise to an isomorphism
    \begin{equation}
      \label{eq:isomj}
      j: Y\xrightarrow{\sim} \mathbb{A}^1, \quad \tau\mapsto j(E_\tau).
    \end{equation}

    \begin{definition}
      Define the surface $X$ to be the product of two modular curves, $$X:=Y\times Y=\{(E,E')\},$$ which is the moduli space of pairs of elliptic curves $(E,E')$. For each positive integer $m$, we define the \emph{modular correspondence} $Z(m)$ over the surface $X$ by $$Z(m):=\{ (E, E',\varphi): \phi: E\xrightarrow{\deg m} E'\}\rightarrow X,$$ parameterizing a pair of elliptic curves $(E,E')$ together with a degree $m$ isogeny $\phi: E\rightarrow E'$. 
    \end{definition}

\begin{figure}[h]
  \centering
  \includegraphics[scale=.7]{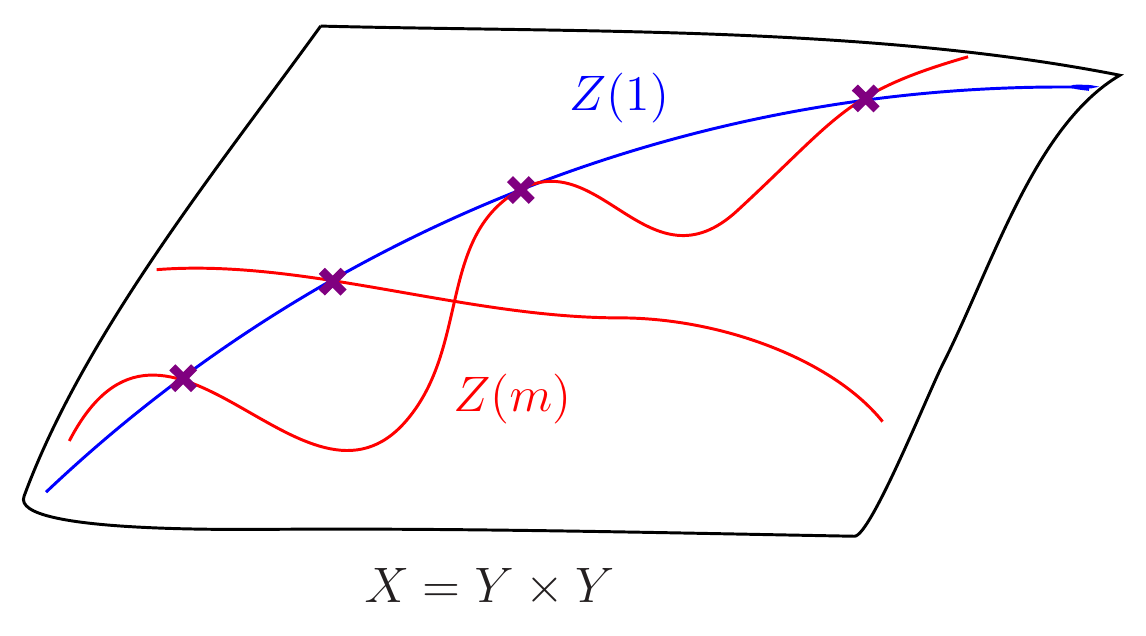}
  \caption{Geometric intersection number}
  \label{fig:hurwitz}
\end{figure}

    The isogeny $\phi$ imposes one nontrivial condition and thus $Z(m)$ defines a divisor on $X$. For example, when $m=1$ the modular correspondence is nothing but the diagonally embedded modular curve $$Z(1)=\Delta (Y)\subseteq X=Y\times Y.$$ Given two divisors $Z(m)$ and $Z(1)$ on the surface $X$, one expects that the intersection $Z(m)\cap Z(1)$ should be 0-dimensional. When this is the case (i.e., when $Z(m)$ and $Z(1)$ intersects \emph{properly}), we obtain a \emph{geometric intersection number} $\langle Z(m),Z(1)\rangle_X$ by counting the number of intersection points weighted by intersection multiplicities (see Figure \ref{fig:hurwitz}). 

Now a curious observation comes: the geometric intersection number $\langle Z(m),Z(1)\rangle_X$ is equal to the LHS of (\ref{eq:hurwitz}),
\begin{equation}
  \label{eq:hurwitz1}
  \langle Z(m),Z(1)\rangle_X=\sum_{dd'=m}\max\{d,d'\}.
\end{equation}
In fact, under the isomorphism (\ref{eq:isomj}) we know that $X$ has a natural compactification $$\overline{X}=\mathbb{P}^1\times \mathbb{P}^1.$$ Since $\overline{Z(1)}=\Delta(\mathbb{P}^1)$, the well-known cohomological equivalence on $\mathbb{P}^1\times \mathbb{P}^1$, $$\Delta(\mathbb{P}^1)\sim \mathbb{P}^1\times \{\text{pt}\}+ \{\text{pt}\}\times \mathbb{P}^1$$  then implies that $$\langle \overline{Z(m)},\overline{Z(1)}\rangle_{\overline{X}}=\langle \overline{Z(m)}, \mathbb{P}^1\times \{\text{pt}\}\rangle_{\overline{X}}+\langle \overline{Z(m)}, \{\text{pt}\}\times\mathbb{P}^1\rangle_{\overline{X}}=2\sum_{d|m}d,$$ where the last equality comes from counting the number of degree $m$ isogenies with a fixed source (resp. target) elliptic curve. The desired identity (\ref{eq:hurwitz1}) then follows from subtracting the contribution at $\{\infty\}=\mathbb{P}^1\setminus \mathbb{A}^1$. 

On the other hand, via the moduli interpretation we have $$Z(m)\cap Z(1)=\{(E,E): E\xrightarrow{\deg m} E\}.$$ When $m$ is a perfect square, we know that  $Z(m)$ contains $Z(1)$ (by considering the multiplication-by-$\sqrt{m}$ isogeny $E\xrightarrow{[\sqrt{m}]} E$), and thus $Z(m)$ and $Z(1)$ do not intersect properly.  However, when $m$ is not a perfect square (as assumed in Theorem \ref{thm:hurwitz}), $E$ has to be an elliptic curve with complex multiplication by an imaginary quadratic order $\mathbb{Z}[\alpha]$, where $\alpha\overline{\alpha}=m$, and thus $Z(m)$ and $Z(1)$ do intersect properly. Using the theory of \emph{complex multiplication}, counting the weighted number of such elliptic curves with complex multiplication exactly gives the sum of Hurwitz class numbers as the RHS of (\ref{eq:hurwitz}),
\begin{equation}
  \label{eq:hurwitz2}
  \langle Z(m),Z(1)\rangle_X=\sum_{t\in \mathbb{Z}, 4m-t^2>0}H(4m-t^2).
\end{equation}
Combining (\ref{eq:hurwitz1}) and (\ref{eq:hurwitz2}) completes our sketch of the geometric proof of the Hurwitz formula (see Gross--Keating \cite{Gross1993} for complete details).

\subsection{Hurwitz formula as a geometric Siegel--Weil formula}\label{sec:geom-sieg-weil} The Hurwitz class numbers appearing in the Hurwitz formula (\ref{eq:hurwitz}) also naturally appear as Fourier coefficients of Siegel Eisenstein series. More precisely, consider the Siegel Eisenstein series $E(\sz,s)$ on $\Sp(4)$ of weight 2, $$E(\sz,s):=\sum_{\left(\begin{smallmatrix}A & B\\ C & D\end{smallmatrix}\right)\in P(\mathbb{Z})\backslash \Sp_4(\mathbb{Z})}\det(C\sz+D)^{-2}\frac{\det\Im(\sz)^{s-1/2}}{|\det(C\sz+D)|^{2s-1}}=\sum_{T\in \Sym_2(\mathbb{Z})_{>0}}E_T(\sz,s)q^T+\cdots.$$ Then we have (up to a normalizing constant) $$H(4m-t^2)=E_T(\tau,s_0=1/2), \quad T=\left(\begin{smallmatrix}m & t/2 \\ t/2 & 1\end{smallmatrix}\right).$$ Notice that the condition $4m-t^2>0$ neatly translates to the condition $T>0$.

We summarize our discussion by the following diagram: 
\begin{equation}
  \label{eq:geometric}
  \xymatrix@R=.5em{\displaystyle\sum_{d d'=m} \max\{ d, d'\}  \ar@{=}[d] & = &  \displaystyle\sum_{t\in \mathbb{Z}, 4m-t^2>0}H(4m-t^2) \ar@{=}[d]\\ \langle Z(m), Z(1)\rangle_{X} & = & \displaystyle\sum_{T=\left(\begin{smallmatrix}m & t/2 \\ t/2 & 1\end{smallmatrix}\right)>0}E_T(\sz,1/2) \\ &&\\ &&\\&&\\
  {\begin{array}[h]{c}
    \text{intersection number of} \\
    \text{2 divisors on } X
  \end{array}}
\ar@{<~>}[uuuu] & = &
  {\begin{array}[h]{c}
   \text{sum of Fourier coefficients of}\\
   \text{a Siegel Eisenstein series on } \Sp(4).
 \end{array}} \ar@{<~>}[uuuu]}
\end{equation}
In this way the Hurwitz formula can be viewed as a \emph{geometric} Siegel--Weil formula for the pair $(G, H)=(\Sp(4), \O(2,2))$, where one replaces the theta integral on $\O(2,2)$ by a \emph{geometric theta series}, i.e., the generating series of geometric intersection numbers of modular correspondences for the surface $X=Y\times Y$,
\begin{center}
  \fbox{\emph{geometric} theta series on $X$ $\longleftrightarrow$ \emph{value} of Siegel Eisenstein series on $\Sp(4)$ at $s_0=1/2$.}
\end{center}
Notice here the natural appearance of the product of modular curves due to the exceptional isomorphism $\SO(2,2)\simeq \SL_2\times \SL_2/\{\pm1\}$. This geometric Siegel--Weil formula further computes a more general geometric intersection number $\langle Z(m), Z(n)\rangle_{X}$ as the sum of $E_T(\sz,1/2)$ with $T=\left(\begin{smallmatrix}m & t/2 \\ t/2 & n\end{smallmatrix}\right)>0$.

\begin{remark}
   The remarkable discovery that generating series involving intersection numbers of cycles are modular originates from the work of Hirzebruch--Zagier \cite{Hirzebruch1976} on Hilbert modular surfaces.  Historically \cite{Hirzebruch1976} was the primary motivation in Kudla's work discussed below  (cf. the introduction of \cite{KM90,Kudla1997}) and also in the work of Gross--Kohnen--Zagier \cite{GKZ87}. See Example \ref{exa:Zm} and Remark \ref{rem:modulairty} (\ref{item:betti}), (\ref{item:chow}). 
\end{remark}

\subsection{Orthogonal Shimura varieties}\label{sec:orth-shim-vari}

Kudla proved a more general \emph{geometric Siegel--Weil formula} by replacing the surface $X=Y\times Y$ by an orthogonal Shimura variety of arbitrary dimension.  Our next goal is to discuss Kudla's formula.

Let $F$ be a totally real number field. Pick a real place $w$ of $F$. Let $V$ be a quadratic space over $F$ of dimension $m$ such that for any place $v|\infty$ of $F$,
\begin{center}
  the $F_v$-quadratic space
  $V_v:= V \otimes_F F_v$ has signature
  $\begin{cases}
   (m-2,2), & \text{if $v=w$}, \\
   (m,0), & \text{if $v\ne w$}.
  \end{cases}$
\end{center}
Let $G=\GSpin(V)$, which sits in an exact sequence $$1\rightarrow \mathbb{G}_m\rightarrow G\rightarrow\SO(V)\rightarrow1.$$ Associated to any open compact subgroup $K\subseteq G(\mathbb{A}_F^\infty)$,  we have a \emph{GSpin Shimura variety} $X=\Sh_G$, which has a smooth canonical model of dimension $m-2$ over the reflex field  $F$ (viewed as a subfield of $\mathbb{C} $ via the embedding induced by the place $w$) and admits complex uniformization $$X(\mathbb{C})=G(F)\backslash[ \mathbb{D}\times G(\mathbb{A}_F^\infty)/K].$$ Here $\mathbb{D}$ is the hermitian symmetric domain of oriented negative 2-planes in $V_w$. Unfolding the definition we may rewrite $X(\mathbb{C})$ as a disjoint union of quotients of $\mathbb{D}$ by congruence subgroups of $G(F)$ (see \cite[\S1]{Kudla2004}). The Shimura variety $X$ is quasi-projective, and is projective when $V$ is anisotropic (e.g., when $F\ne \mathbb{Q}$, by the signature condition).

\begin{remark}
One technical reason that one prefers to work with the Shimura variety $X$ associated to $G$ (instead of $\SO(V)$) is that it is of \emph{Hodge type} (instead of \emph{abelian type}), and admits an embedding into a Siegel modular variety (the moduli space of polarized abelian varieties) of larger dimension.  
\end{remark}

\begin{example}(cf. \cite{Kudla2004,Howard2014,FH00}) Consider $F=\mathbb{Q}$. Via accidental isomorphisms between $G=\GSpin(m-2,2)$ and classical groups of symplectic type in low ranks, the Shimura varieties $X$ recover many classical modular varieties in low dimensions (see Table \ref{tab:GSpin},  where $D$ is a quaternion algebra over $\mathbb{Q}$ and $E=\mathbb{Q}(\sqrt{d})$ is a real quadratic field). When $\dim X\le 19$, $X$ is also closely related to the moduli space of polarized K3 surfaces and has proved to be useful for studying the arithmetic of K3 surfaces.
    \begin{table}[h]
    \begin{tabular}[h]{|c|c|c|}
      $\dim X$ & $G=\GSpin(m-2,2)$ & $X$\\\hline
      1 & $\GL_2$ or $D^\times$ & modular/Shimura curve\\
      2 & $\GL_2\times_{\mathbb{G}_m} \GL_2$ or $\GL_{2,E}^{\det \in \mathbb{Q}^\times}$ & product of modular curves or Hilbert modular surfaces\\
      3 & $\mathrm{GSp}_4$  & Siegel 3-fold = moduli of abelian surfaces\\
      4 & $\mathrm{GU}(2,2)$ (up to center) & moduli of abelian 4-folds with complex multiplication\\
      6 & $\mathrm{GSp}(4, D)$ & moduli of abelian 8-folds with quaternion multiplication\\
    \end{tabular}
    \centering\caption{Examples of GSpin Shimura varieties.  }\label{tab:GSpin}
  \end{table}
\end{example}

\subsection{Kudla's generating series of special cycles and the modularity conjecture}\label{sec:kudl-gener-seri}

The Shimura variety $X$ is equipped with special divisors $Z(m)\rightarrow X$ generalizing  in the case $X=Y\times Y$. Via an embedding into a Siegel modular variety, $Z(m)$ parameterizes certain polarized abelian varieties together with a \emph{special endomorphism} of degree $m$ (see \cite[\S5]{MadapusiPera2016}).

\begin{example}\label{exa:Zm}
 The special divisor $Z(m)$ recovers
  \begin{itemize}
  \item Heegner points on modular curves and Shimura curves  when $\dim X=1$ (cf. \cite[Appendix]{Kudla2004}),
  \item modular correspondences $Z(m)$ on $X=Y\times Y$ considered in \S\ref{sec:geometric-proof} when $\dim X=2$ and $G$ is split,
  \item Hirzebruch--Zagier cycles \cite{Hirzebruch1976} on Hilbert modular surfaces when $\dim X=2$ and $G$ is nonsplit.
  \end{itemize}
\end{example}

More generally, for any $y\in V$ with $(y,y)>0$, its orthogonal complement $V_y\subseteq V$ has rank $n-1$. The embedding  $G_y:=\GSpin(V_y)\hookrightarrow G=\GSpin(V)$ defines a Shimura subvariety of codimension 1 $$\Sh_{G_y}\rightarrow X=\Sh_{G}.$$ For any $x\in V(\mathbb{A}_F^\infty)$ with $(x,x)\in F_{>0}$, there exists $y\in V$ and $g\in G(\mathbb{A}_F^\infty)$ such that $y=gx$. Define the \emph{special divisor} $$Z(x)\rightarrow X$$ to be the $g$-translate of $\Sh_{G_y}$. For any $\mathbf{x}=(x_1,\ldots, x_\m)\in V(\mathbb{A}_F^\infty)^\m$ with $(\mathbf{x},\mathbf{x})\in \Sym_{\m}(F)_{>0}$, define the \emph{special cycle} (of codimension $\m$) $$Z(\mathbf{x})=Z(x_1)\cap \cdots \cap Z(x_\m)\rightarrow X.$$ Here $\cap$ denotes the fiber product over $X$. More generally, for a $K$-invariant Schwartz function $\varphi\in \sS(V(\mathbb{A}_F^\infty)^\m)^K$ and $T\in \Sym_{\m}(F)_{>0}$,  define the \emph{weighted special cycle} $$Z_\varphi(T)=\sum_{\mathbf{x}\in K\backslash V(\mathbb{A}_F^\infty)^\m\atop (\mathbf{x},\mathbf{x})=T}\varphi(\mathbf{x}) Z(\mathbf{x})\in \Ch^\m(X)_\mathbb{C}:=\Ch^\m(X)\otimes \mathbb{C}.$$ Here $\Ch^\m(X)$ is the Chow group of algebraic cycles of codimension $\m$ on $X$ (up to rational equivalence). With extra care, we can also define $Z_\varphi(T)\in \Ch^\m(X)_\mathbb{C}$ for any $T\in \Sym_{\m}(F)_{\ge0}$ (see \cite{Kudla2004,YZZ09}).

\begin{definition}
Define \emph{Kudla's generating series of special cycles}
\begin{equation}
  \label{eq:kudlagen}
  Z_\varphi(\sz)=\sum_{T\in\Sym_\m(F)_{\ge0}}Z_\varphi(T)q^T,
\end{equation}
 as a formal sum valued in $\Ch^\m(X)_\mathbb{C}$, where $$\sz\in \mathcal{H}_n=\{\sx+i\sy: \sx\in\Sym_n(F_\infty),\ \sy\in\Sym_n(F_\infty)_{>0}\},\quad q^T:=\prod_{v|\infty}e^{2\pi i \tr T\sz_v}.$$
\end{definition}

\begin{remark}
  Analogous constructions of special divisors and Kudla's generating series also apply to Shimura varieties of unitary type associated to hermitian spaces with signature $(m-1,1)$ at one archimedean place and signature $(m,0)$ at all other archimedean places (see \cite{Liu2011}). These Shimura varieties of orthogonal/unitary type can be naturally viewed as Shimura varieties associated to totally definite incoherent quadratic/hermitian spaces (\cite{Zha19,Gro20}), see \S\ref{sec:unit-shim-vari}. In the unitary case, we obtain a generating series of the form
  \begin{equation}
    \label{eq:hermgen}
    Z_\varphi(\sz)=\sum_{T\in\Herm_\m(F)_{\ge0}}Z_\varphi(T)q^T,
  \end{equation}
where we replace positive semi-definite symmetric matrices $\Sym_n(F)_{\ge0}$ by positive semi-definite hermitian matrices $\Herm_n(F)_{\ge0}$, and replace Siegel's half space by the hermitian half space $$\sz\in\mathcal{H}_n:=\{\sx+i\sy: \sx\in\Herm_n(F_\infty),\ \sy\in\Herm_n(F_\infty)_{>0}\}.$$
\end{remark}

 We may view $Z_\varphi(\sz)$ as a geometric theta series, now valued in Chow groups for cycles of arbitrary codimension $\m$. The analogy to Siegel's theta series \eqref{eq:siegeltheta} and the theta integral (\ref{eq:thetaintegral}) leads to the following Kudla's modularity conjecture.

  \begin{conjecture}[Kudla's modularity] \label{conj:modularity}
    The formal generating series $Z_\varphi(\sz)$ converges absolutely and defines a modular form on $\mathcal{H}_n$ of weight $m/2$ valued in $\Ch^\m(X)_\mathbb{C}$. 
  \end{conjecture}
  \begin{remark}\quad\label{rem:modulairty}
    \begin{altenumerate}
    \item\label{item:betti} The analogous modularity in Betti cohomology, i.e., the modularity of the generating series valued in $H^{2\m}(X(\mathbb{C}), \mathbb{C})$ defined by the image of $Z_\varphi(T)$  under the cycle class map $$\Ch^\m(X)\rightarrow H^{2\m}(X(\mathbb{C}), \mathbb{Z}),$$ is known by the classical work of Kudla--Millson \cite{KM90}. The special case of special divisors on Hilbert modular surfaces dates back to Hirzebruch--Zagier \cite{Hirzebruch1976} (see also Funke--Millson \cite{FM14}).
    \item\label{item:chow} Kudla's modularity conjecture was originally formulated for orthogonal Shimura varieties over $\mathbb{Q}$ (\cite{Kudla1997,Kudla2004}). In this case, Borcherds \cite{Bor99} proved the conjecture for the divisor case $m=1$ (the special case of Heegner points on modular curves dates back to the classical work of Gross--Kohnen--Zagier \cite{GKZ87}). Zhang \cite{Zha09} proved the modularity for general $\m$ assuming the absolute convergence of the series. Bruinier--Westerholt-Raum \cite{BW15} proved the desired convergence and hence established Kudla's modularity conjecture for orthogonal Shimura varieties over $\mathbb{Q}$. More recently, Bruinier--Zemel \cite{BZ19} has extended the modularity to toroidal compactifications of orthogonal Shimura varieties when $\m=1$.
    \item For orthogonal Shimura varieties over totally real fields, Yuan--Zhang--Zhang \cite{YZZ09} proved the modularity for $\m=1$ (see also Bruinier \cite{Bru12} for a different proof) and reduce the $\m>1$ case to the convergence.
    \item Conjecture~\ref{conj:modularity} in the unitary case was formulated by Liu \cite{Liu2011}, who also proved the case $\m=1$ and reduce the $\m>1$ case to the convergence. Recently Xia \cite{Xia21} proved the desired convergence when $E=\mathbb{Q}(\sqrt{-d})$ for $d=1,2,3,7,11$ (in the notation of \S\ref{sec:unit-shim-vari}), and thus established Conjecture~\ref{conj:modularity} in these cases.
    \item\label{item:arithmetic} Kudla \cite[Problem 4]{Kudla2004} also proposed the modularity problem in the \emph{arithmetic Chow group} $\aCh^\m(\mathcal{X})$ of a suitable (compactified) integral model $\mathcal{X}$ of $X$ (see \cite{GS90,BGKK07} and also \cite{Sou92}). The problem seeks to define canonically an  explicit arithmetic generating series  $\widehat{\mathcal{Z}}_\varphi(\sz)$ valued in $\aCh^\m(\mathcal{X})_\mathbb{C}$ which lifts $Z_\varphi(\sz)$ under the restriction map $$\aCh^\m(\mathcal{X})\rightarrow\Ch^\m(X),$$ and such that $\widehat{\mathcal{Z}}_\varphi(\sz)$ is modular. When $\m=1$, this \emph{arithmetic modularity} was proved  by Howard--Madapusi Pera \cite{HMP20} (orthogonal groups over $\mathbb{Q}$) and Bruinier--Howard--Kudla--Rapoport--Yang \cite{Bruinier2017} (unitary groups over $\mathbb{Q}$). Several low dimensional cases were also proved:
      \begin{itemize}
      \item Shimura/modular curves (Kudla--Rapoport--Yang \cite{Kudla2006}, Sankaran \cite{San14} and Du--Yang \cite{DY19}),
      \item Hilbert modular surfaces (Bruinier--Burgos Gil--K\"uhn \cite{BBGK07}),
      \item Product of modular curves (Berndt--K\"uhn \cite{BK12,BK12a}).
      \end{itemize}
 We also mention the arithmetic modularity of the difference of two arithmetic theta series by Ehlen--Sankaran \cite{ES18} for $\m=1$ (unitary groups over $\mathbb{Q}$), the almost arithmetic modularity by Mihatsch--Zhang \cite[Theorem 4.3]{MZ21} for $\m=1$ (unitary groups over totally real fields $F\ne \mathbb{Q}$), the arithmetic modularity of Fourier--Jacobi coefficients for general $\m$ by Sankaran \cite{San20} (anisotropic orthogonal groups), and several striking recent works  involving applications of arithmetic modularity \cite{Andreatta2018, Zhang2019, SSTT19}.
    \item Finally, we mention several recent works on the modularity conjecture for more general classes of orthogonal and unitary Shimura varieties (indefinite at more than one archimedean places) by Rosu--Yott \cite{RY20}, Kudla \cite{Kud19} and Maeda \cite{Mae19,Mae21}.
    \end{altenumerate}
  \end{remark}

\subsection{Kudla's geometric Siegel--Weil formula}

In the special case $n=\dim X$, the generating series $Z_\varphi(\sz)$ in \eqref{eq:kudlagen} is valued in $\Ch^{\dim X}(X)_\mathbb{C}$ (i.e., the Chow group of 0-cycles). When $X$ is projective, composing with the degree map $$\deg:\Ch^{\dim X}(X)_\mathbb{C}\rightarrow \mathbb{C},$$ we obtain a generating series $\deg Z_\varphi(\sz)$ valued in $\mathbb{C}$. Its terms encode geometric intersection numbers between special divisors on $X$ generalizing the case $X=Y\times Y$ considered in \S\ref{sec:geom-sieg-weil}. Kudla  proved the following remarkable geometric version of the Siegel--Weil formula (by analogy with Theorem \ref{thm:SWF} specialized to the case $n=\dim X=m-2$ and so $s_0=1/2$).

\begin{theorem}[Kudla's geometric Siegel--Weil formula {\cite[Corollary 10.5]{Kudla1997}}]
  Assume that $X$ is projective (i.e., $V$ is anisotropic). Take $n=\dim X$. The for any $\varphi\in \sS(V(\mathbb{A}_F^\infty)^\m)^K$ the following identity holds (up to a nonzero constant depending only on choices of measures)
  \begin{equation}
    \label{eq:kudla}
    \deg Z_\varphi(\sz)\doteq E(\sz, 1/2,\varphi \otimes \varphi_\infty).
  \end{equation}
  Here $\varphi_\infty\in \sS(V(F_\infty)^\m)$ is a certain Schwartz function constructed from the Kudla--Millson Schwartz form (\cite{Kudla1986}).
\end{theorem}

Thus Kudla's geometric Siegel--Weil formula is a precise identity of the form \begin{center}
\fbox{\emph{geometric} theta series  on $X$ $\longleftrightarrow$ \emph{value} of Siegel Eisenstein series on $\Sp(2n)$ at $s_0=1/2$.}
\end{center}

\begin{remark}
Kudla in fact proved a geometric Siegel--Weil formula for the generating series of special cycles of all dimension, i.e., without assuming $n=\dim X$. Here and in the next section we focus on the case of 0-cycles and refer to Kudla's excellent surveys \cite{Kudla2004, Kud02,Kud02a} for the general case.
\end{remark}

\section{Arithmetic Siegel--Weil formula}

In this section we discuss an arithmetic version of the Siegel--Weil formula. Parallel to the previous section, we will use an example in the case $(G,H)=(\Sp(6),\O(2,2))$ considered by Gross--Keating to motivate the more general case.

\subsection{Gross--Keating formula}\label{sec:gross-keat-form}

\begin{figure}[h]
  \centering
  \includegraphics[scale=.7]{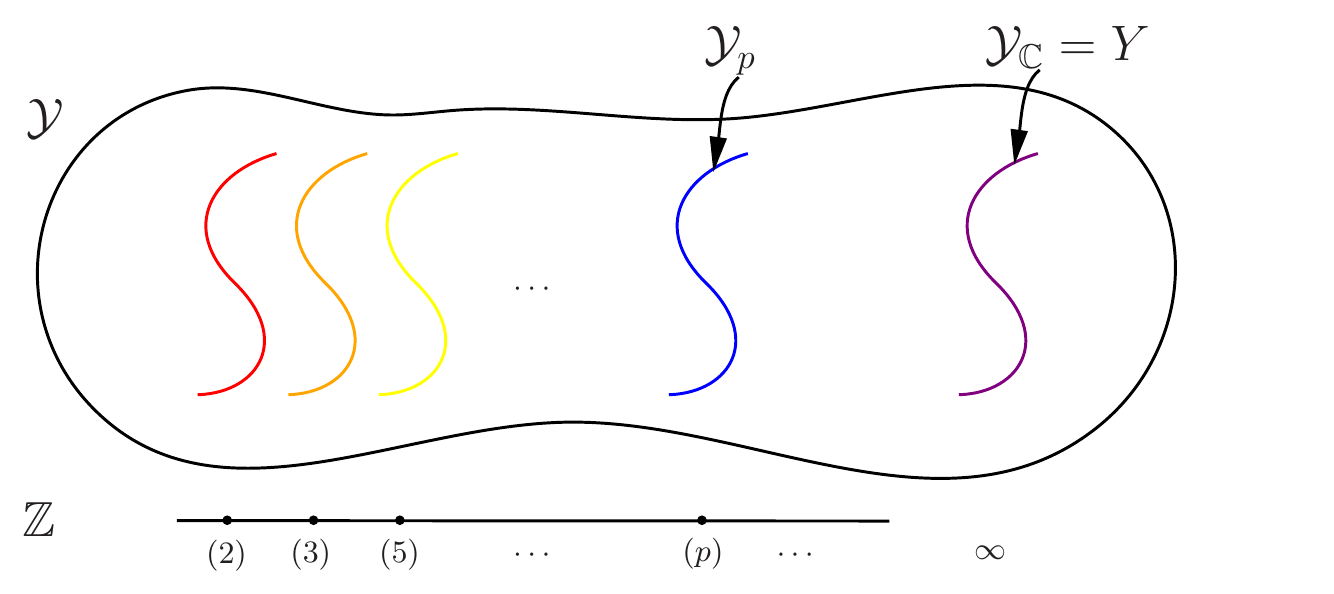}
  \caption{Arithmetic surface}
  \label{fig:curve}
\end{figure}

Gross--Keating took the geometric point of view of the Hurwitz class number relation and found a remarkable generalization for \emph{arithmetic} intersection numbers. As the moduli space of elliptic curves, the modular curve $Y$ has a canonical integral model $\mathcal{Y}$ over $\Spec \mathbb{Z}$ such that $\mathcal{Y}(\mathbb{C})=Y$. The integral model $\mathcal{Y}$ is an \emph{arithmetic surface} fibered over the arithmetic curve $\Spec \mathbb{Z}$, and its fiber above $p$ is a smooth curve in characteristic $p$ (Figure \ref{fig:curve}).

Analogously, the surface $X$ have an canonical integral model $\mathcal{X}$ over $\Spec \mathbb{Z}$, which is an \emph{arithmetic threefold}. The modular correspondence $Z(m)\rightarrow X$ naturally extends to a divisor $\mathcal{Z}(m)\rightarrow \mathcal{X}$. Now on the arithmetic threefold $\mathcal{X}$, we need 3 (instead of 2) divisors so that the intersection has expected dimension 0. Define the \emph{arithmetic intersection number} by $$\langle \mathcal{Z}(m_1), \mathcal{Z}(m_2), \mathcal{Z}(m_3)\rangle_{\mathcal{X}}:=\sum_p \langle \mathcal{Z}(m_1),\mathcal{Z}(m_2), \mathcal{Z}(m_3)\rangle_{\mathcal{X}_p}\cdot \log p,$$ where $\langle \mathcal{Z}(m_1),\mathcal{Z}(m_2), \mathcal{Z}(m_3)\rangle_{\mathcal{X}_p}$ encodes the intersection number supported in the fiber $\mathcal{X}_p$. This definition makes sense when these three divisors intersect properly (or more generally when their intersection is supported in finitely many fibers $\mathcal{X}_p$).

Is there a formula for this arithmetic intersection number analogous to (\ref{eq:geometric})? The index set for $T$ should be $\Sym_3(\mathbb{Z})$ to put the three diagonal entries $m_1,m_2,m_3$ in $T$, and thus we need to look at the Siegel Eisenstein series $E(\sz,0)$ on $\Sp(6)$ of weight 2 (instead of on $\Sp(4)$). Then the relevant special point in the Siegel--Weil formula is $s_0=0$, the central point. However, it is too naive to expect $\langle \mathcal{Z}(m_1), \mathcal{Z}(m_2), \mathcal{Z}(m_3)\rangle_{\mathcal{X}}$ to be equal to $E(\sz,0)$: the arithmetic intersection number involves a $\mathbb{Q}$-linear combination of $\log p$'s and hence is no longer a rational number like $E_T(\sz,1/2)$ in (\ref{eq:geometric}). Moreover, the Eisenstein series $E(\sz,s)$ turns out to have an odd functional equation at the center $s=0$, and hence $E(\sz,0)=0$ automatically! This automatic vanishing suggests that it is interesting to look at its first \emph{derivative} $E'(\sz,0)$ at $s=0$.

After these two appropriate modifications: replacing $\Sp(4)$ by $\Sp(6)$ and replacing the value at $s_0=1/2$ by the central derivative at $s_0=0$,  it turns out that  we do have the following remarkable formula relating the arithmetic intersection numbers to $E'(\sz,0)$ (see a nice exposition of the proof in \cite{VGW+07}).

\begin{theorem}[Gross--Keating \cite{Gross1993}, Gross--Kudla--Zagier] Assume there is no positive definite binary quadratic form representing $m_1,m_2,m_3$ simultaneously. Then (up to an explicit constant)
  \begin{equation}
    \label{eq:grosskeating}
    \xymatrix@R=4em{\langle \mathcal{Z}(m_1), \mathcal{Z}(m_2), \mathcal{Z}(m_3)\rangle_{\mathcal{X}} & \doteq &  \displaystyle\sum_{T=\left(\begin{smallmatrix}m_1 & * & *\\ * & m_2 & * \\ * & * &m_3 \end{smallmatrix}\right)>0} E_T'(\sz,0)\\
  {\begin{array}[h]{c}
    \text{arithmetic intersection number}\\
    \text{of 3 divisors on } \mathcal{X}
  \end{array}} \ar@{<~>}[u]
 & &
{\begin{array}[h]{c}
  \text{sum of Fourier coefficients of the derivative}\\
\text{ of a Siegel Eisenstein series on } \Sp(6)
\end{array}}  \ar@{<~>}[u]
} 
  \end{equation}
\end{theorem}

\begin{remark}\label{rem:properint}
  The assumption on $m_1,m_2,m_3$ is analogous to the assumption that $m$ is not a perfect square in Theorem \ref{thm:hurwitz}, which guarantees that the three divisors intersect properly.
\end{remark}

\subsection{Arithmetic Siegel--Weil formula} 

Parallel to (\ref{eq:geometric}),  the  Gross--Keating formula can be viewed as an \emph{arithmetic} Siegel--Weil formula for the pair $(G, H)=(\Sp(6), \O(2,2))$, where one replaces the theta integral on $\O(2,2)$ by a generating series of arithmetic intersection numbers on the arithmetic threefold $\mathcal{X}$,
\begin{center}
\fbox{\emph{arithmetic} theta series on $\mathcal{X}$ $\longleftrightarrow$ \emph{central derivative} of Siegel Eisenstein series on $\Sp(6)$.}
\end{center}

Of course there is nothing stopping us from considering the higher dimensional case. In fact Kudla (\cite{Kudla1997a}) and Kudla--Rapoport (\cite{Kudla1999a,Kudla2000a,Kudla2014}) proposed vast conjectural generalizations of the Gross--Keating formula by
\begin{enumerate}
\item Taking $\mathcal{X}$ to be a suitable regular integral model $\mathcal{X}$ of the Shimura variety $X$ of orthogonal/unitary type considered in \S\ref{sec:orth-shim-vari}--\ref{sec:kudl-gener-seri} associated to quadratic/hermitian spaces. 
\item Defining suitable integral models $\mathcal{Z}(m)\rightarrow \mathcal{X}$  of the special divisors $Z(m)\rightarrow X$, and more generally integral models $\mathcal{Z}(T)$ for special divisors $Z(T)$.
\end{enumerate}

\begin{remark}
  Such a regular integral model $\mathcal{X}$ is known to exists for certain level subgroups $K$ defining $X$ (cf. \cite[\S14]{LZ} in the unitary case). To avoid technicalities we will be vague about the choice of level structures and integral models, but see Remark \ref{sec:local-arithm-sieg-1} (\ref{item:ramified}), (\ref{item:parahoric}). 
\end{remark}

Now take $n=\dim \mathcal{X}$ ( $=\dim X+1$) so that the intersection of $n$ special divisors has expected dimension 0. We have a natural decomposition $$\mathcal{Z}(m_1)\cap \cdots \cap\mathcal{Z}(m_n)=\bigsqcup_{T=\left(\begin{smallmatrix}m_1 & * & *\\ * & \cdots & * \\ * & * &m_n \end{smallmatrix}\right)} \mathcal{Z}(T)$$ index by symmetric/hermitian matrices $T$ with diagonal entries $m_1,\ldots, m_n$. Here $\cap$ denotes the fiber product over $\mathcal{X}$. When $T>0$, it turns out that $\mathcal{Z}(T)$ is supported in finitely many fibers $\mathcal{Z}(T)_p$, and we have a well-defined $T$-part of the arithmetic intersection number $$\langle \mathcal{Z}(m_1), \cdots, \mathcal{Z}(m_n)\rangle_T:=\sum_p \langle \mathcal{Z}(m_1), \cdots, \mathcal{Z}(m_n)\rangle_{\mathcal{Z}(T)_p}\cdot\log p.$$   
Now we are ready to state the conjecture on an arithmetic version of the Siegel--Weil formula (by analogy with Theorem~\ref{thm:SWF} in the orthogonal case specialized to $n=m+1$ and so $s_0=0$ is the central point), which is known as the \emph{Kudla--Rapoport conjecture} in the unitary case ({\cite[Conjecture 11.10]{Kudla2014}}).

\begin{conjecture}[Arithmetic Sigel--Weil formula, nonsingular part]\label{conj:ASW} Take $n=\dim \mathcal{X}$. Then for any $T\in \Sym_n(F)_{>0}$ (resp. $T\in \Herm_n(F)_{>0}$) in the orthogonal (resp. unitary) case with diagonal entries $m_1,\ldots, m_n$, the following identity holds (up to a nonzero constant depending only on choices of measures),
  \begin{equation}
    \label{eq:ASW}
      \xymatrix@R=1em{\langle \mathcal{Z}(m_1),\cdots, \mathcal{Z}(m_n)\rangle_T
 & \stackrel{?}{=}& E_T'(\sz,0)\\ {\begin{array}[h]{c}
    \text{arithmetic intersection number}\\
    \text{of $n$ divisors on } \mathcal{X}
  \end{array}} \ar@{<~>}[u]
 & &
{\begin{array}[h]{c}
  \text{Fourier coefficient of the central derivative of a}\\
\text{  Siegel Eisenstein series on } \Sp(2n) \text{(resp. $\UU(n,n)$)}
\end{array}}  \ar@{<~>}[u]
}.
  \end{equation}
\end{conjecture}

\begin{remark}\label{rem:improper}
  In general the $n$ special divisors do not intersect properly, and a more sophisticated definition of the arithmetic intersection numbers is needed (cf. Definition \ref{def:intL}). In particular, with the correct definition the conjecture works even for improper intersections.
\end{remark}

Thus the arithmetic Siegel--Weil formula a precise conjectural identity of the form
\begin{center}
\fbox{\emph{arithmetic} theta series on $\mathcal{X}$ $\longleftrightarrow$ \emph{central derivative} of Siegel Eisenstein series on $\Sp(2n)$ ($\UU(n,n)$).}  
\end{center}

Now we can state one of the main results of \cite{LZ} (see \cite[Theorem 1.3.1]{LZ} for more precise technical assumptions).

\begin{theorem}[with Zhang \cite{LZ}]\label{thm:arithm-sieg-weil}
Conjecture \ref{conj:ASW} holds for arbitrary $n$ in the unitary case.
\end{theorem}

Our recent work with Zhang \cite{LZ22} has also established a slightly weaker semi-global (at a good odd prime $p$) version of  Conjecture \ref{conj:ASW} in the orthogonal case. We will discuss some key ideas of the proof in \S\ref{sec:local-arithm-sieg}.

\begin{remark}\quad
  \begin{altenumerate}
  \item   When $n=3$ and $G$ is split (orthogonal case), the arithmetic Siegel--Weil formula recovers the Gross--Keating formula, and also treats the case of improper intersections (see Remarks \ref{rem:properint}, \ref{rem:improper}).
  \item  Theorem \ref{thm:arithm-sieg-weil} was previously proved when $n=3$ by Terstiege \cite{Terstiege2011, Terstiege2013} and when the intersection is 0-dimensional (\cite{Gross1993,Kudla1999a,Kudla2000a,Kudla2011, Bruinier2018}).
  \item There is also an \emph{archimedean} part of the arithmetic Siegel--Weil formula, relating archimedean arithmetic intersection numbers with the nonsingular but indefinite Fourier coefficients of $E'(\sz,0)$. These Fourier coefficients are nonholomorphic, unlike the positive definite Fourier coefficients in Conjecture \ref{conj:ASW}. This archimedean arithmetic Siegel--Weil formula was proved by Liu \cite{Liu2011} (unitary case), and Garcia--Sankaran \cite{Garcia2019} in full generality (see also Bruinier--Yang \cite{Bruinier2018} for an alternative proof in the orthogonal case).
  \item   Kudla conjectured that there should also be a \emph{singular} part of the arithmetic Siegel--Weil formula, relating the singular Fourier coefficients of $E'(\sz,0)$ to certain arithmetic intersection numbers. However the singular part is more difficult to prove, or even to formulate precisely, cf. \cite[Problem 6]{Kudla2004}. As a special case, the constant term of the arithmetic Siegel--Weil formula should roughly relate the arithmetic volume of $\mathcal{X}_n$ to logarithmic derivatives of Dirichlet $L$-functions. Such an explicit arithmetic volume formula was proved by H\"ormann \cite{H14} (orthogonal case) and Bruinier--Howard \cite{BH21} (unitary case), though a precise comparison with the constant term of $E'(\sz,0)$ is yet to be formulated and established.
  \item Ideally, putting all singular/nonsingular and archimedean/nonarchimedean parts together, one should arrive at a full arithmetic Siegel--Weil formula of the form in complete analogy to (\ref{eq:kudla}), $$\widehat{\deg}\ \widehat{\mathcal{Z}}_\varphi(\sz)\stackrel{?}{=} E'(\sz,0,\varphi \otimes \varphi_\infty).$$   Here $\widehat{\mathcal{Z}}_\varphi(\sz)$ is the conjectural arithmetic theta series in Remark \ref{rem:modulairty} (\ref{item:arithmetic}), $\widehat{\deg}$ is the arithmetic degree map $$\widehat{\deg}: \aCh^{\dim \mathcal{X}}(\mathcal{X})_\mathbb{C}\rightarrow \mathbb{C},$$ and $\varphi_\infty$ is the standard Gaussian function on the totally positive definite space over $F_\infty$.  The full arithmetic Siegel--Weil formula was established by Kudla, Rapoport and Yang (\cite{Kudla1999, Kudla1997a, Kudla2000, Kudla2006}) for $n=1,2$ (orthogonal case) in great generality.  However, it remains an open problem to formulate such a precise full arithmetic Siegel--Weil formula in higher dimension.
  \item Recently Feng--Yun--Zhang \cite{FYZ21} proved a higher Siegel--Weil formula over function fields for unitary groups, which relates nonsingular coefficients of the $r$-th derivative of Siegel Eisenstein series and intersection numbers of special cycles on moduli spaces of Drinfeld shtukas with $r$ legs.  The case $r=0$ (resp. $r=1$) can be viewed as an analogue of the Siegel--Weil formula (resp. the arithmetic Siegel--Weil formula). Over function fields, the possibility of relating higher derivatives of analytic objects to intersection numbers was first discovered by Yun--Zhang \cite{YZ17,YZ19} in the context of the higher Gross--Zagier formula.  Over number fields, however, no analogue of such a higher Siegel--Weil formula (resp. higher Gross--Zagier formula) is currently known when $r>1$. Feng--Yun--Zhang  \cite{FYZ21a} also defined higher theta series over function fields (including all singular terms) and conjectured their modularity. Previously, an arithmetic Siegel--Weil formula over function fields was proved by Wei \cite{Wei19} for special cycles on moduli spaces of Drinfeld modules of rank 2 with complex multiplication (analogue of the special case $n=r=1$).
  \end{altenumerate}
\end{remark}

\section{Local arithmetic Siegel--Weil formula}\label{sec:local-arithm-sieg}

In order to prove the arithmetic Siegel--Weil formula (Conjecture \ref{conj:ASW}), one first notices that it can be reduced to a local identity:
\begin{enumerate}
\item Geometric side (LHS): The arithmetic intersection numbers corresponding to nonsingular matrices $T$ can be computed as a sum indexed by primes of $\mathbb{Z}$ (or the ring of integers $O_F$ in general). The local term at a finite prime $p$ can be further reduced to an arithmetic intersection on a Rapoport-Zink space $\mathcal{N}_n$, which is a local analogue of Shimura varieties over $\mathbb{Z}_p$ (or a completion of $O_F$ in general), via the theory of $p$-adic uniformization of Shimura varieties (\cite{RZ96}).
\item Analytic side (RHS): The nonsingular Fourier coefficients $E_T(\sz,s)$ has a product expansion indexed by primes of $\mathbb{Z}$, and thus the derivative $E_T'(\sz,0)$  can also be written as a sum  indexed by the primes of $\mathbb{Z}$. The term indexed by a finite prime $p$ can be further reduced to the derivative of the local representation density of quadratic/hermitian forms over $\mathbb{Z}_p$.
\end{enumerate}

This reduction step is illustrated in the following diagram:
 $$\xymatrix@R=4em{
{\begin{array}[h]{c}
    \text{arithmetic intersection number}\\
    \text{of $n$ special divisors on } \mathcal{X}
  \end{array}}
\ar@{<~>}[d]
 & \stackrel{?}{=}&
{\begin{array}[h]{c}
  \text{sum of Fourier coefficients of central derivative}\\
\text{ of a Siegel Eisenstein series on }  \Sp(2n)  \text{ or } \UU(n,n)
\end{array}}
\ar@{<~>}[d]\\
{\begin{array}[h]{c}
    \text{arithmetic intersection number}\\
    \text{of $n$ special divisors on } \mathcal{N}_n
  \end{array}}
& \stackrel{?}{=}&
{\begin{array}[h]{c}
  \text{central derivative of the local density}\\
\text{ of a hermitian or quadratic form in $n$-variables}
\end{array}}
}$$

The conjectural local identity on the bottom is known as the \emph{local arithmetic Siegel--Weil formula}. The local arithmetic Siegel--Weil formula has been recently proved in our work with Zhang \cite{LZ} (resp.  \cite{LZ22}) in the unitary (resp. orthogonal) case. Next we will make this local conjecture more precise in the unitary case. In the unitary case,  this reduction step was made precise by Kudla--Rapoport (\cite{Kudla2011}, \cite{Kudla2014}) and the local conjecture is also known as the \emph{local Kudla--Rapoport conjecture}. 

\subsection{Geometric side}\label{sec:geometric-side} Let $p$ be an odd prime. Let $F_0$ be a finite extension of $\mathbb{Q}_p$ with residue field $k=\mathbb{F}_q$ and uniformer $\varpi$. Let $F$ be the unramified quadratic extension of $F_0$ (e.g., $F/F_0=\mathbb{Q}_{p^2}/\mathbb{Q}_p$). Associated to this datum we have the\emph{ unitary Rapoport-Zink space $\mathcal{N}_n$}:

\begin{itemize}
\item It is a formal scheme over $\Spf_{\OFb}$ of relative dimension $n-1$, parameterizing deformations (up to quasi-isogenies) of a fixed $O_F$-hermitian formal $p$-divisible group $\mathbb{X}_n/\bar k$ of relative height $2n$, dimension $n$ and signature $(1,n-1)$. Here $\OFb$ is the completion of the maximal unramified extension of $O_F$.
\item The space of \emph{special homomorphisms} $$\mathbb{V}_n=\Hom_{O_F}(\mathbb{X}_1, \mathbb{X}_n)\otimes_{O_F} F$$ has a structure of a (non-split) $F$-hermitian space of dimension $n$, coming from the principal polarization on $\mathbb{X}_n$.
\item The unitary group $\UU(\mathbb{V}_n)$ naturally acts on $\mathcal{N}_n$ via the action on $\mathbb{X}_n$.
\item Each vector $x\in \mathbb{V}_n$ gives rise to a \emph{special divisor} or \emph{Kudla--Rapoport (KR) divisor} $\mathcal{Z}(x)\subseteq \mathcal{N}_n$, defined to be the locus where the homomorphism $x$ deforms. This is the local analogue of the special divisor considered in Conjecture \ref{conj:ASW}.
\end{itemize}

The Rapoport--Zink space $\mathcal{N}_n$ is formally smooth over $\Spf \OFb$, but its geometric structure is rather complicated. For example, $\mathcal{N}_n$ is highly non-reduced: the reduced subscheme $\mathcal{N}_n^\mathrm{red}$ has dimension $\lfloor \frac{n-1}{2}\rfloor$, near the middle dimension of $\mathcal{N}_n$. The structure of $\mathcal{N}^\mathrm{red}_n$ was studied by Vollaard--Wedhorn (\cite{Vollaard2011}), and they showed that $\mathcal{N}_n^\mathrm{red}$ has a nice stratification into smooth varieties, known as the \emph{Bruhat--Tits stratification}. Each closed stratum of the Bruhat--Tits stratification is isomorphic to a smooth projective variety $\DL_i$ of dimension $i$ over $\bar k$ ($0\le i\le\lfloor \frac{n-1}{2}\rfloor$), and the incidence relation between the closed strata resembles the combinatorial structure of the Bruhat--Tits building for quasi-split unitary groups  over $F_0$. Here each $\DL_i$ is a \emph{generalized Deligne--Lusztig variety} associated to the unitary group $\UU(2i+1)$ over $k$.

\begin{example} Take $n=1$. Then $\mathcal{N}_1\cong\Spf \OFb$, and $\mathcal{N}_1^\mathrm{red}=\{\text{pt}\}$.
  \end{example}

\begin{example}
  Take $n=3$. Then $\mathcal{N}_3$ has relative dimension 2 over $\Spf \OFb$, while the reduced subscheme $\mathcal{N}_3^\mathrm{red}$ has dimension 1 (Figure \ref{fig:RZ}). In this case only two types of Deligne--Lusztig variety show up:
  \begin{enumerate}
  \item   $\DL_0=\{\text{pt}\}$, a single point.
  \item  $\DL_1=\{x^{q+1}+y^{q+1}+z^{q+1}=0\}\subseteq \mathbb{P}^2$, the Fermat curve of degree $q+1$.
  \end{enumerate}
  And $\mathcal{N}_3^\mathrm{red}$ is an infinite tree, where
  \begin{enumerate}
  \item the number of $\DL_1$ containing a given $\DL_0$ is exactly $q+1$.
  \item the number of $\DL_0$ contained in a given $\DL_1$ is exactly $q^3+1$.
  \end{enumerate}
  \begin{figure}[h]
    \centering
    \includegraphics[scale=1.1]{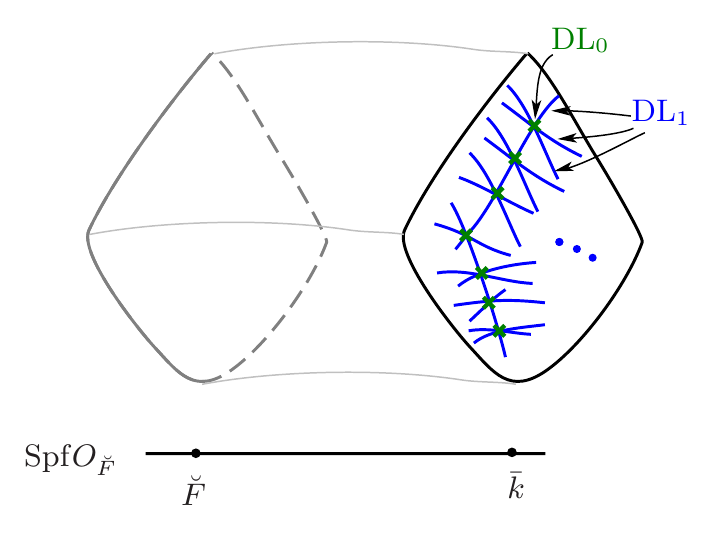}
    \caption{$n=3$}
    \label{fig:RZ}
  \end{figure}
\end{example}

\begin{definition}\label{def:intL}
  Let $L\subseteq \mathbb{V}_n$ be an $O_F$-lattice of rank $n$. Let $x_1,\ldots, x_n$ be an $O_F$-basis of $L$. Define the \emph{special cycle} or \emph{Kudla--Rapoport (KR) cycle} $$\mathcal{Z}(L):=\mathcal{Z}(x_1)\cap\cdots\cap \mathcal{Z}(x_n)\subseteq \mathcal{N}_n.$$ Define the \emph{arithmetic intersection number} $$\Int(L) :=\chi(\mathcal{O}_{\mathcal{Z}(x_1)} \otimes^\mathbb{L}\cdots \otimes^\mathbb{L}\mathcal{O}_{\mathcal{Z}(x_n)} ),$$ where $\chi$ denotes the Euler--Poincar\'e characteristic, $\mathcal{O}_{\mathcal{Z}(x_i)}$ denotes the structure sheaf of the Kudla--Rapoport divisor $\mathcal{Z}(x_i)$, and $\otimes^\mathbb{L}$ denotes the derived tensor product of coherent sheaves on $\mathcal{N}_n$.    It is known (by Terstiege \cite{Terstiege2013} as extended in \cite[Corollary 2.8.2]{LZ}, or by Howard \cite{Howard2018}) that $\Int(L)$ is independent of the choice of the basis $x_1,\ldots, x_n$ and hence is a well-defined invariant of $L$ itself, justifying the notation.
\end{definition}

\begin{remark}
When the intersection $\mathcal{Z}(L)$ is 0-dimensional, we have $$\mathcal{O}_{\mathcal{Z}(x_1)}\otimes^\mathbb{L}\cdots \otimes^\mathbb{L}\mathcal{O}_{\mathcal{Z}(x_n)}=\mathcal{O}_{\mathcal{Z}(x_1)}\otimes\cdots \otimes\mathcal{O}_{\mathcal{Z}(x_n)}=\mathcal{O}_{\mathcal{Z}(L)}$$ and thus $\Int(L)$ is nothing but the $\OFb$-length of $\mathcal{O}_{\mathcal{Z}(L)}$ (= the sum of intersection multiplicities at all points). Even though  $\mathcal{Z}(L)$ is the intersection of $n$ divisors in a $n$-dimensional formal scheme, in general $\mathcal{Z}(T)$ may be fail to have the expected dimension 0 due to improper intersection. In this case, the derived intersection is needed so that the intersection number is well-behaved.
\end{remark}

\begin{example}\label{exa:IntL}
  Take $n=3$ and $L\cong\langle \varpi\rangle^3$ (the hermitian form with respect to an $O_F$-basis of $L$ is $\diag\{\varpi,\varpi,\varpi\}$). In this case the intersection $\mathcal{Z}(L)$ is not 0-dimensional, and in fact $$\mathcal{Z}(L)\cong \DL_1=\{x^{q+1}+y^{q+1}+z^{q+1}=0\}\subseteq \mathbb{P}^2.$$ The arithmetic intersection number turns out to be its topological Euler characteristic $$\Int(L)=\chi(\DL_1)=2-2g=2-q(q-1)=2+q-q^2.$$
\end{example}

\subsection{Analytic side}

\begin{definition}
Let $L, M$ be two hermitian $O_F$-lattices of rank $n, m$ respectively. Let $S=\Rep(M, L)$ be the $O_{F_0}$-scheme such that for any $O_{F_0}$-algebra $R$, $$S(R)=\mathrm{Herm}(L \otimes_{F_0}R, M \otimes_{F_0}R),$$ where $\mathrm{Herm}$ denotes the set of hermitian module homomorphisms. The \emph{local density} of representations of $M$ by $L$ is defined to be $$\Den(M,L):=\lim_{N\rightarrow +\infty}\frac{\# S(O_{F_0}/\varpi^N)}{q^{N\cdot\dim S_{F_0}}}.$$ It gives a quantitative measure of ``how many different ways'' one can embed $L$ into $M$ as a hermitian submodule.   
\end{definition}

\begin{example}
  Consider $L=M=\langle 1\rangle^n$, the rank $n$ lattice with hermitian form given by the identity matrix. Then one can compute that $$\Den(\langle 1\rangle^n, \langle 1\rangle^n)=\prod_{i=1}^{n}(1-(-q)^{-i}).$$  One can recognize it as the number $\frac{G(k)}{q^{\dim G}}$, where $G$ is the unitary group $\U(\langle 1\rangle^n)$ defined over $O_{F_0}$, or in fancier language, as the local $L$-factor of the Gross motive \cite{Gro97} of the quasi-split unitary group in $n$ variables.
\end{example}

The local density $\Den(M, L)$ has nice compatibility when replacing $M$ by $M \oplus \langle 1\rangle^k$ for $k\ge0$. More precisely, it is known (\cite[Theorem II]{Hironaka1998}) that $\Den(\iden^{n+k}, L)$ is a polynomial in $(-q)^{-k}$ with $\mathbb{Q}$-coefficients.

\begin{example}[{\cite[p.677]{Kudla2011}}]
  \begin{equation}
  \label{eq: iden}
\Den(\iden^{n+k}, \iden^n)=\prod_{i=1}^n(1-(-q)^{-i}X)\bigg|_{X= (-q)^{-k}}.
\end{equation}
\end{example}

\begin{definition}
  Define the \emph{normalized Siegel series} $\Den(X,L)\in \mathbb{Z}[X]$ such that $$\Den((-q)^{-k},L)=\frac{\Den(\iden^{n+k}, L)}{\Den(\iden^{n+k}, \iden^n)}.$$
\end{definition}

These polynomials $\Den(X,L)$ plays an important role in computing the Fourier coefficients of Siegel Eisenstein series, and many works are devoted to proving more explicit formulas for them (see e.g., \cite{Kit83, Katsurada1999, Hironaka1998, Hironaka2012,IK16,CY}). The local Siegel series satisfies a functional equation (\cite[Theorem 5.3]{Hironaka2012})
\begin{equation}
  \label{eq:functionalequation}
  \Den(X,L)=(-X)^{\val(L)}\cdot \Den\left(\frac{1}{X},L\right).
\end{equation} Here $\val(L):=\val(\det(L))\in \mathbb{Z}$ is the \emph{valuation} of $L$. This may be viewed as a local analogue of the functional equation of Eisenstein series.

\begin{definition}
If $\val(L)$ is odd (equivalently, $L \otimes_{O_F}F$ is a non-split hermitian space), then $\Den(1,L)=0$ by the functional equation (\ref{eq:functionalequation}). In this case, define the \emph{central derivative of the local density} by $$\pDen(L)\coloneqq-\frac{\rd}{\rd X}\bigg|_{X=1}\Den(X,L).$$ 
\end{definition}

\subsection{Local arithmetic Siegel--Weil formula} Now we are ready to state the local arithmetic Siegel--Weil formula in the unitary case, originally conjectured by Kudla--Rapoport \cite[Conjecture 1.3]{Kudla2011}.

\begin{theorem}[{Local arithmetic Siegel--Weil formula, with Zhang \cite[Theorem 1.2.1]{LZ}}]\label{thm:KR}
    Let $L\subseteq \mathbb{V}_n$ be an $O_F$-lattice of full rank $n$. Then $$\Int(L)=\pDen(L).$$
  \end{theorem}

\begin{remark}\quad\label{sec:local-arithm-sieg-1}
  \begin{altenumerate}
  \item   The theorem was proved by Kudla--Rapoport  \cite{Kudla2011} for $n=2$ and Terstiege \cite{Terstiege2013} for $n=3$.
 \item   The local arithmetic Siegel--Weil formula is proved for the orthogonal case in our work with Zhang  \cite{LZ22} for arbitrary $n$. The case $n=3$ was previously proved by Gross--Keating \cite{Gross1993} and Terstiege \cite{Terstiege2011}.
 \item\label{item:ramified}   The local arithmetic Siegel--Weil formula is proved when the quadratic extension $F/F_0$ is ramified for exotic smooth models in our work with Liu \cite{LL2020} for arbitrary even $n$, and for Kr\"amer models\footnote{During the refereeing process of this article, He--Shi--Yang \cite{HSY21} has formulated a conjectural local arithmetic Siegel--Weil formula for Kr\"amer models for arbitrary $n$ and proved it for the case $n=3$. The case for arbitrary $n$ has been proved in our work with He--Shi--Yang \cite{HLSY22}.} by Shi and He--Shi--Yang \cite{Shi20, HSY20} for $n=2$. 
 \item\label{item:parahoric}   It is more difficult to prove or formulate the local arithmetic Siegel--Weil formula in the presence of more general level structures (even when the quadratic extension $F/F_0$ is unramified). In the unitary case \cite{LZ} formulates and proves a local arithmetic Siegel--Weil formula when the level is the parahoric subgroup given by the stabilizer of an almost self-dual lattice (the case $n=2$ was previously proved by Sankaran \cite{Sankaran2017}).  Recently Cho \cite{Cho20} has proposed a general formulation for all minuscule parahoric levels in the unitary case.
  \end{altenumerate}
\end{remark}

\begin{example}
  Take $n=3$ and $L\cong\langle \varpi\rangle^3$. Specializing the formula of Cho--Yamauchi \cite{CY}  (extended to the unitary case in \cite[Theorem 3.5.1]{LZ}) gives
  \begin{align*}
    \Den(X,L)&=(1-X)(1+qX)(1-q^2X)+(q^3+1)(1-X)X^2\\
    & = 1-(1-q+q^2)X+(1-q+q^2)X^2-X^3
  \end{align*}
It satisfies the functional equation $$\Den(X,L)=-X^3\cdot \Den\left(\frac{1}{X},L\right).$$ It is easy to compute
\begin{align*}
  \pDen(L)&=-\frac{\rd}{\rd X}\bigg|_{X=1}\Den(X,L)\\
  &= (1-q+q^2)-2(1-q+q^2)+3\\
  &= 2+q-q^2.
\end{align*}
So combining Example \ref{exa:IntL} we obtain $\Int(L)=\pDen(L)$ in this case! It is miraculous that the purely analytic quantity $\pDen(L)$ secretly knows about the Euler characteristic $\chi(\DL_1)=\Int(L)$ of the Deligne--Lusztig curve $\DL_1$.
\end{example}

\subsection{Strategy of the proof: uncertainty principle}
The previously known special cases of Theorem \ref{thm:KR} were proved via explicit computation of both the geometric and analytic sides. Explicit computation seems infeasible for the general case. The proof in \cite{LZ} instead proceeds via induction on $n$ using the \emph{uncertainty principle}, a standard tool from local harmonic analysis. Even for $n=2,3$, this proof is different from the previous proofs.

More precisely, for a fixed $O_F$-lattice $L^\flat\subseteq \mathbb{V}:=\mathbb{V}_n$ of rank $n-1$, consider functions on $x\in \mathbb{V}\setminus L^\flat_F$, $$\mathrm{Int}_{L^\flat}(x):= \mathrm{Int}(L^\flat+\langle x\rangle),\quad \partial\mathrm{Den}_{L^\flat}(x):= \partial\mathrm{Den}(L^\flat+\langle x\rangle).$$ Then it remains to show the equality of the two functions $$\mathrm{Int}_{L^\flat}=\partial\mathrm{Den}_{L^\flat}.$$ By the definition of $\Int(L)$ and $\pDen(L)$, both functions are easily seen to vanish when $x$ is non-integral, i.e., $\mathrm{val}(x)<0$. Here $\mathrm{val}(x)$ denotes the valuation of the norm of $x$. By utilizing the inductive structure of the Rapoport--Zink spaces and local densities, it is not hard to see that if $x\perp L^\flat$ with $\mathrm{val}(x)=0$, then $$\mathrm{Int}_{L^\flat}(x)=\mathrm{Int}(L^\flat),\quad \partial\mathrm{Den}_{L^\flat}(x)=\partial\mathrm{Den}(L^\flat)$$ for the lattice $L^\flat\subseteq \mathbb{V}_{n-1}\cong \langle x\rangle^\perp_F$ of full rank $n-1$. Thus by induction on $n$, the difference function $\phi=\mathrm{Int}_{L^\flat}-\partial\mathrm{Den}_{L^\flat}$ already vanishes on a large subset $$\{x \in\mathbb{V}: x\perp L^\flat, \mathrm{val}(x)\le0\}.$$ We would like to deduce that $\phi$ indeed vanishes identically. To this end, we apply the following Uncertainty Principle.

\begin{proposition}[Uncertainty Principle, {\cite[Proposition 8.1.6]{LZ}}]
  Let $\phi\in \sS(\mathbb{V})$ be a Schwartz function on $\mathbb{V}$. If both $\phi$ and its Fourier transform $\hat\phi$ vanish on $\{x\in \mathbb{V}: \mathrm{val}(x)\le0\}$. Then $\phi=0$.
\end{proposition}

In other words, $\phi,\hat\phi$ cannot simultaneously have ``small support'' unless $\phi=0$. Applying the Uncertainty Principle to the difference function $\phi$, then we can finish the proof as long as we get a good control over the support of $\hat \phi$. However, both functions $\mathrm{Int}_{L^\flat},\partial\mathrm{Den}_{L^\flat}$ have \emph{singularities} along the hyperplane $L^\flat_F:=L^\flat \otimes_{O_F}F\subseteq \mathbb{V}$. Intuitively, if $x$ is ``closer'' to $L^\flat$, then $\mathcal{Z}(x)$ and $\mathcal{Z}(L^\flat)$ intersect more ``improperly'', which results in the blow-up of $\Int_{L^\flat}(x)$ along $L^\flat_F$. These singularities cause trouble in computing the Fourier transforms or even in showing that $\phi\in \sS(\mathbb{V})$.

\subsection{Strategy of the proof: decomposition and local modularity}

\begin{figure}[h]
  \centering
  \includegraphics{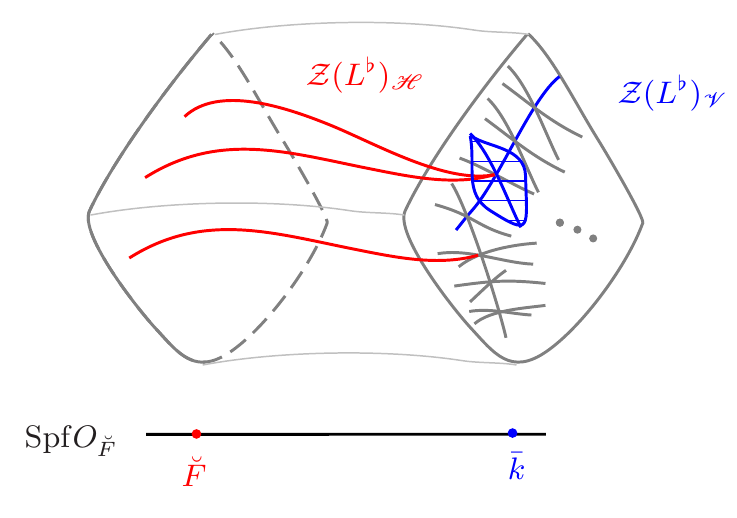}
  \caption{Decomposition}
  \label{fig:decomp}
\end{figure}

To overcome this difficulty, we isolate the singularities by decomposing $$\mathrm{Int}_{L^\flat}=\mathrm{Int}_{L^\flat,\sH}+\mathrm{Int}_{L^\flat,\sV},\quad \partial\mathrm{Den}_{L^\flat}=\partial\mathrm{Den}_{L^\flat,\sH}+\partial\mathrm{Den}_{L^\flat,\sV}$$ into ``horizontal'' and ``vertical'' parts. Here on the geometric side $\mathrm{Int}_{L^\flat,\sH}$ (resp. $\mathrm{Int}_{L^\flat,\sV}$) is the contribution from the horizontal (resp. vertical part) of the KR cycles, illustrated by the red (resp. blue) part in Figure \ref{fig:decomp}. One can hope to understand the horizontal part explicitly using deformation theory, and the vertical part using algebraic geometry over the residue field $\bar k$.

By the uncertainty principle and the induction on $n$, it remains to prove the following.
  \begin{theorem}[Key Theorem]\label{thm:Key}\quad
    \begin{enumerate}[label=(KR\arabic*)]
    \item $\Int_{L^\flat,\sH}=\pDen_{L^\flat,\sH}$.
    \item $\Int_{L^\flat,\sV}\in \sS(\mathbb{V})$\quad and\quad  $\widehat{\Int_{L^\flat,\sV}}=- \Int_{L^\flat,\sV}.$
    \item $\pDen_{L^\flat,\sV}\in \sS(\mathbb{V})$\quad and \quad  $\widehat{\pDen_{L^\flat,\sV}}=- \pDen_{L^\flat,\sV}.$
    \end{enumerate}
  \end{theorem}

In other words,  the horizontal part $\partial\mathrm{Den}_{L^\flat,\sH}$ matches $\mathrm{Int}_{L^\flat,\sH}$, and subtracting the horizontal parts removes the singularities along $L^\flat_F$ so that vertical parts are indeed in $\sS(\mathbb{V})$. The extra invariance under Fourier transform
\begin{equation}
  \label{eq:localmodular}
  \widehat{\mathrm{Int}}_{L^\flat,\sV}=- \mathrm{Int}_{L^\flat,\sV},\quad \widehat{\partial\mathrm{Den}}_{L^\flat,\sV}=-\partial\mathrm{Den}_{L^\flat,\sV}.
\end{equation}
 can be thought of as a \emph{local modularity}, by analogy with the global modularity of arithmetic generating series (such as in Bruinier--Howard--Kudla--Rapoport--Yang \cite{Bruinier2017} discussed in Remark \ref{rem:modulairty} (\ref{item:arithmetic})) encoding an extra global $\SL_2$-symmetry.

\subsection{Ingredients of the Key Theorem}

Some key ingredients of the proof of the Key Theorem~\ref{thm:Key} include:
\begin{enumerate}[label=(KR\arabic*)]
\item\label{item:KR1} We describe explicitly the horizontal part of KR cycles  in terms of Gross's quasi-canonical liftings \cite{Gross1986a}, using the work of Tate, Grothendieck--Messing and Breuil on the deformation theory of $p$-divisible groups.
\item\label{item:KR2} On the geometric side we show (\ref{eq:localmodular}) by reducing to the case of intersection with Deligne--Lusztig curves $\DL_1$. This reduction requires the Bruhat--Tits stratification of $\mathcal{N}^\mathrm{red}$ into the Deligne--Lusztig varieties $\DL_i$, as discussed in \S\ref{sec:geometric-side}, and the Tate conjecture \cite{Tate1994} for these Deligne--Lusztig varieties. We prove the latter by reducing to a cohomological computation of Lusztig \cite{Lusztig1976/77}.
\item\label{item:KR3} On the analytic side we show (\ref{eq:localmodular}) using  Cho--Yamauchi's explicit formula \cite{CY} for $\pDen(L)$  in terms of weighted lattice counting, and reduce to a (rather subtle) lattice theoretic problem. (In fact we only show directly something weaker than (\ref{eq:localmodular}) which is enough to imply Theorem \ref{thm:KR}, and we then deduce (\ref{eq:localmodular}) a posteriori).
\end{enumerate}

\section{Arithmetic inner product formula}

In this last section we discuss an application of the arithmetic Siegel--Weil formula to the Beilinson--Bloch conjecture and the arithmetic inner product formula.

\subsection{Birch and Swinnerton-Dyer conjecture}

One long-standing problem in number theory is the determination of the rational points $E(\mathbb{Q})$ for an elliptic curve $E: y^2=x^3 +a x+b$ defined over $\mathbb{Q}$. The celebrated \emph{Birch and Swinnerton-Dyer (BSD) conjecture} predicts a deep link between $E(\mathbb{Q})$ and its $L$-function $L(E,s)$. Define the \emph{algebraic rank} $$r_\mathrm{alg}(E):=\rank E(\mathbb{Q})$$ to be the rank of the finitely generated abelian group $E(\mathbb{Q})$. Define the \emph{analytic rank}  $$r_\mathrm{an}(E):=\ord_{s=1} L(E,s)$$ to be the order of vanishing of $L(E,s)$ at the central point $s=1$. The BSD conjecture predicts the rank equality between these two notions of ranks of seemingly different nature,
\begin{equation}
  \label{eq:rank}
  r_\mathrm{an}(E)\stackrel{?}{=}r_\mathrm{alg}(E),
\end{equation}
 and a refined BSD formula
 \begin{equation}
   \label{eq:BSDformula}
   \frac{L^{(r)}(E,1)}{r!}\stackrel{?}=\frac{\Omega(E) R(E) \prod_p c_p(E)\cdot |\Sha(E)|}{|E(\mathbb{Q})_\mathrm{tor}|^2}
 \end{equation}
 for the Taylor expansion of $L(E/\mathbb{Q},s)$ at $s=1$ (here $r=r_\mathrm{an}(E)$) in terms of various important arithmetic invariants of $E$. Among these invariants are the order of the (mysterious) \emph{Tate--Shafarevich group} $\Sha(E)$, and the \emph{regulator} $R(E):=\det (\langle P_i, P_j\rangle_{\text{NT}})_{r\times r},$  where  $$\langle \ , \ \rangle_\text{NT}: E(\mathbb{Q})\times E(\mathbb{Q})\rightarrow \mathbb{R}$$ is the\emph{ N\'eron--Tate height pairing} and $\{P_i\}$ is a basis of the free part of $E(\mathbb{Q})$.

 The BSD conjecture is still widely open in general, but much progress has been made in the rank 0 or 1 case. The seminal work of Gross--Zagier \cite{Gross1986}, Kolyvagin \cite{Kol90} proved the implications
 \begin{equation}
   \label{eq:GZK}
   r_\mathrm{an}(E)=0\Rightarrow r_\mathrm{alg}(E)=0,\quad r_\mathrm{an}(E)=1\Rightarrow r_\mathrm{alg}(E)=1,
 \end{equation}
 confirming (\ref{eq:rank}) when $r_\mathrm{an}(E)\le1$. Due to the work of many people, many cases of (\ref{eq:BSDformula}) are also known when $r_\mathrm{an}(E)\le1$.

 The key to relate $r_\mathrm{an}(E)$ and $r_\mathrm{alg}(E)$ is the\emph{ Gross--Zagier formula}
 \begin{equation}
   \label{eq:GZ}
   L'(E,1)\doteq \langle P, P\rangle_\mathrm{NT},
 \end{equation}
 (up to an explicit nonzero constant, including the period of $E$ and other rational factors) relating the first derivative of the $L$-function at $s=1$ and the N\'eron--Tate height of certain rational points $P$ on $E$ known as \emph{Heegner points} (cf. Example \ref{exa:heegner}). It gives one crucial implication in (\ref{eq:GZK}),
 \begin{equation}
   \label{eq:GZimply}
   r_\mathrm{an}(E)=1\Longrightarrow r_\mathrm{alg}(E)\ge1.
 \end{equation}
 The tools of Heegner points and $L$-functions, linked via the Gross--Zagier formula, are indispensable in studying the arithmetic of elliptic curves. We refer to Zhang \cite{Zhang2014a} for an excellent recent survey on Heegner points and the Birch--Swinnerton-Dyer conjecture (see also Gross \cite{Gro04} and Darmon \cite{Dar04}).

 \subsection{Beilinson--Bloch conjecture}

 In the 1980s, Beilinson (\cite[Conjecture 5.9]{Beuilinson1987}) and Bloch (\cite[Recurring Fantasy]{Bloch1984},\cite[Conjecture]{Bloch1984a}) proposed vast generalizations of the BSD conjecture to higher dimensional varieties.

 Let $X$ be a smooth projective variety of dimension $n$ over a number field $F$. For $1\le \i\le n$, denote by $\Ch^\i(X)$ the Chow group of algebraic cycles of codimension $\i$ on $X$ defined over $F$ (up to rational equivalence), and $\Ch^\i(X)^0\subseteq \Ch^\i(X)$ the subgroup of geometrically cohomologically trivial cycles. Denote by $L(H^{2\i-1}(X), s)$ the motivic $L$-function associated to the $(2\i-1)$-th \'etale cohomology $H^{2\i-1}(X_{\overline{F}}, \mathbb{Q}_\ell)$. Then the \emph{Beilinson--Bloch (BB) conjecture} predicts the equality between analytic and algebraic ranks
\begin{equation}
  \label{eq:BB}
\ord_{s=\i}L(H^{2\i-1}(X), s)= \rank \Ch^\i(X)^0,
\end{equation}
and a refined formula $$L^{(r)}(H^{2\i-1}(X),\i)\doteq \det \langle \ , \ \rangle_{\mathrm{BB}}$$ for the leading coefficient at $s=\i$ in terms of the determinant of the Beilinson--Bloch height pairing\footnote{When $\i>1$, the Beilinson--Bloch height pairing is only defined assuming certain conjectures on algebraic cycles on $X$ (see \cite[Conjectures 2.2.1 and 2.2.3]{Beuilinson1987}). This important technical issue is addressed in \cite{LL2020, LL2021} so that the right-hand-side of (\ref{eq:AIPF}) in Theorem \ref{thm:AIPF} can be defined unconditionally, but we will intentionally ignore it for the purpose of this article.} $$\langle\ ,\  \rangle_\mathrm{BB}: \Ch^\i(X)^0\times \Ch^{n+1-\i}(X)^0\rightarrow \mathbb{R}.$$

\begin{example}
When $X/F=E/\mathbb{Q}$ and $\i=1$ we recover the BSD conjecture as
\begin{center}
  $\Ch^1(E)^0\simeq E(\mathbb{Q})$, $L(H^1(E), s)=L(E,s)$ and
  $\langle\ ,\ \rangle_\mathrm{BB}=-\langle\ ,\ \rangle_\mathrm{NT}$.
\end{center}
\end{example}

 The BB conjecture is even more elusive than the BSD conjecture: in general we do not know that $\Ch^\i(X)^0$ is finitely generated, nor do we know that $L(H^{2\i-1}(X), s)$ can be analytically continued to the central point $s=\i$,  so neither side of (\ref{eq:BB}) is well-defined! This may be more an exciting challenge than disappointment for mathematicians --- after all  we were in a similar circumstance when BSD conjecture was formulated in the 1960s: we knew neither the analytic continuation of $L(E,s)$ (except when $E$ has complex multiplication) nor the finiteness of $\Sha(E)$ in order to make sense of either side of (\ref{eq:BSDformula}).

A good testing ground for the BB conjecture is by taking $X$ to be Shimura varieties. By Langland's philosophy, their $L$-functions can be computed in terms of automorphic $L$-functions, so the analytic rank in (\ref{eq:BB}) becomes accessible. Even though we do not known if $\Ch^\i(X)^0$ is finitely generated, it still makes sense to test if it is nontrivial. It is thus tempting to relate special cycles on Shimura varieties to  automorphic $L$-functions, and in particular, to generalize the Gross--Zagier formula (\ref{eq:GZ}) to higher dimensional Shimura varieties and prove the analogue of (\ref{eq:GZimply})  towards the BB conjecture
\begin{equation}
  \label{eq:nonvanishing}
\ord_{s=\i}L(H^{2\i-1}(X), s)=1\Longrightarrow \rank\Ch^\i(X)^0\ge 1.  
\end{equation} Here we use the notation $\rank\Ch^\i(X)^0\ge 1$ to stand for the non-triviality of $\Ch^\i(X)^0_\mathbb{Q}$.

\subsection{Arithmetic inner product formula}\label{sec:arithm-inner-prod}

As discussed in \S\ref{sec:kudl-gener-seri}, Shimura varieties $X$ of type $\GSpin(n-1,2)$ and $\U(n-1,1)$ admit a rich supply of special cycles, recovering Heegner points in the case of modular curves. As explained in \cite[III]{Kudla2004}, when $n$ is even, the arithmetic Siegel--Weil formula together with the doubling method of Piatetski-Shapiro--Rallis \cite{PR84} has important application to the \emph{arithmetic inner product formula} of the form $$L'(1/2, \pi)\doteq\langle \Theta_\varphi(\phi), \Theta_\varphi(\phi)\rangle_\mathrm{BB}.$$ Here
\begin{itemize}
\item  $\pi$ is a certain cuspidal automorphic representation on $\Mp(n)$ or $\UU(n)$.
\item $L'(1/2,\pi)$ is the central \emph{derivative} the standard $L$-function $L(s,\pi)$ of $\pi$ (cf. \cite{Yam14}) with global root number $\varepsilon(\pi)=-1$.
\item $\Theta_\varphi(\phi)$ is an algebraic cycle  on $X$ of codimension $n/2$, constructed from holomorphic forms $\phi\in\pi$ using the method of \emph{arithmetic theta lifting} (see Definition \ref{def:arithm-theta-lift}).
\end{itemize}

This arithmetic inner product formula can be viewed as an arithmetic analogue of the \emph{Rallis inner product formula} (see \cite{GQT14}) of the form
\begin{equation}
  \label{eq:rallis}
  L(1/2, \pi)\doteq\langle \theta_\varphi(\phi),\theta_\varphi(\phi)\rangle_\mathrm{Pet}
\end{equation}
  which relates the  central \emph{value} $L(1/2,\pi)$ when $\varepsilon(\pi)=+1$ and  the Petersson inner product of the (usual) theta lift $\theta_\varphi(\phi)$ discussed in Remark \ref{rem:thetalifting}. It can also be viewed as a higher dimensional generalization of the Gross--Zagier formula (the case $n=2$). Thus the arithmetic Siegel--Weil formula is intimately linked to the Birch--Swinnerton-Dyer conjecture, and more generally the Beilinson--Bloch conjecture for higher dimensional Shimura varieties.

The conjectural arithmetic inner product formula was formulated by Kudla \cite{Kudla1997a} using the Gillet--Soul\'e height and in more generality by Liu \cite{Liu2011} using the Beilinson--Bloch height (in the unitary case). In the unitary case the arithmetic inner product formula has been recently proved under mild local assumptions and Kudla's modularity conjecture in our works with Liu \cite{LL2020, LL2021}, which has been applied to prove the first unconditional theorem for the Beilinson--Bloch conjecture for higher dimensional Shimura varieties. Our remaining goal is to explain some details about the main results of \cite{LL2020, LL2021}.

\subsection{Arithmetic theta lifting on modular curves} 
To motivate the construction of arithmetic theta lifting, let us first consider an example of Heegner points on elliptic curves. 

\begin{example}[cf. \cite{Zag85}]\label{exa:heegner}
  Consider the curve $37a1$ in Cremona's table, $$E=37a1: y^2 + y = x^3 - x,$$ It is the rank one optimal curve over $\mathbb{Q}$ of smallest conductor ($N=37$). It is isomorphic to the modular curve $X_0^+(N)=X_0(N)/\langle w_{N}\rangle$ for $N=37$ ($N$ is the smallest such that $X_0^+(N)$ has positive genus). The Mordell-Weil group $E(\mathbb{Q})\cong \mathbb{Z}$, generated by $P=(0,0)$.  Let $$f_E=\sum_{n\ge1}a_nq^n=q - 2 q^2 - 3 q^3 + 2 q^4 - 2 q^5 + 6 q^6 - q^7 + 6 q^9 + 4 q^{10} - 5 q^{11} - 6 q^{12}  +\cdots \in S_2(N)$$ be the weight 2 newform of level $N=37$ associated to $E$. It gives a modular parametrization $$\varphi_E: \mathcal{H}\rightarrow E=\mathbb{C}/\Lambda,\quad \varphi_E(\tau)=\sum_{n\ge1}\frac{a_nq^n}{n}=q-q^2-q^3+\frac{1}{2}q^4-\frac{2}{5}q^5+q^6-\frac{1}{7}q^7+\cdots $$ such that $\varphi_E^*(\omega_E)=2\pi if_E(\tau)d\tau$. The Shimura--Waldspurger correspondence, which can be viewed as an instance of the theta correspondence for the pair $(G,H)=(\mathrm{Mp}(2),\O(3))$,  $$\theta: S_{3/2}^+(4N) \rightarrow S_2(N),$$ gives a weight 3/2 newform in Kohnen's plus space $S_{3/2}^+(4N)$, $$g_E=\sum_{d\ge1}c_dq^d=-q^3-q^4+q^7-q^{11}+q^{12}+2q^{16}+3q^{27}+\cdots -6q^{67}+\cdots$$ such that $\theta(g_E)=f_E$.

Let $d$ be a positive integer such that $-d\equiv0,1\pmod{4}$. If $N=37$ splits in $K=\mathbb{Q}(\sqrt{-d})$, then one can construct a rational point $P_d\in E(\mathbb{Q})$ using Heegner points on $X_0(N)$ associated to quadratic orders of discriminant $-d$. For example, when $-d$ is a fundamental discriminant, we have a Heegner point $y_K\in X_0(H_K)$ defined over the Hilbert class field $H_K$ of $K$, and $P_d=\tr_{H_K/K}\varphi_E(y_K-\infty)$.  The point $P_d$ may depend on the choice of $d$, even when $E(\mathbb{Q})\cong\mathbb{Z}$. In Table \ref{tab:heegner},  we compute a list of Heegner points $P_d\in E(\mathbb{Q})$ for small $d$'s, and also the integers multiples $n_d$ such that $P_d=n_d\cdot P$ for the generator $P=(0,0)$.
  \begin{table*}[h]
    \centering
    \begin{tabular}[h]{*{11}{|>{$}c<{$}}|}
      d & 3 & 4 & 7 & 11 & 12 & 16 & 27 & \cdots &67& \cdots\\\hline
      P_d & (0,-1) & (0,-1) & (0,0) & (0,-1) & (0,0) & (1,0) & (-1,-1) &\cdots &(6,-15) &\cdots\\\hline      
      n_d & -1 & -1 & 1 & -1 & 1 & 2 & 3 &\cdots &-6 &\cdots\\      
    \end{tabular}
    \caption{Heegner points on $E=X_0^+(37)$}\label{tab:heegner}
  \end{table*}
  
From Table \ref{tab:heegner}, we observe the miraculous coincidence that the integer $n_d$ exactly matches the coefficient $c_d$ of $g_E$! In other words, the generating series $$Z(\tau)=\sum_{d\ge1}P_d q^d=g_E\cdot P\in S_{3/2}^+(4\cdot 37) \otimes E(\mathbb{Q})$$ is a \emph{modular form} valued in $E(\mathbb{Q}) \otimes \mathbb{C}$. This allows us to define the \emph{arithmetic theta lift} by taking the Petersson inner product of $Z(\tau)$ with $g_E$, $$\Theta(g_E):=\langle Z(\tau), g_E\rangle_{\mathrm{Pet}}=\langle g_E,g_E\rangle_\mathrm{Pet }\cdot P \in E(\mathbb{Q})\otimes \mathbb{C},$$ which is now a canonical element in $E(\mathbb{Q}) \otimes \mathbb{C}$ (i.e., no need to choose any particular $d$). The arithmetic inner product formula in this case asserts the identity
  \begin{equation}
    \label{eq:heegnerAIPF}
    L'(E,1)\doteq\frac{\langle \Theta(g_E), \Theta(g_E)\rangle_\mathrm{NT}}{\langle g_E, g_E\rangle_\mathrm{Pet}}.
  \end{equation}
 In fact, in this case we have $\frac{\langle \Theta(g_E), \Theta(g_E)\rangle_\mathrm{NT}}{\langle g_E, g_E\rangle_\mathrm{Pet}}=\langle g_E, g_E\rangle_\mathrm{Pet}\cdot \langle P, P\rangle_\mathrm{NT}$ and we can compute explicitly
  \begin{align*}
    \langle g_E, g_E\rangle_\mathrm{Pet}&=0.7146356107\cdots=\frac{3\omega^+(E)}{4\pi},\\
    \langle P, P\rangle_\mathrm{NT}&= 0.0511114082\cdots,\\
    L'(E,1)&= 0.3059997738    \cdots.
  \end{align*}
Here $\omega^+(E)$ is the real period of $E$. Therefore $$\frac{L'(E,1)}{\langle g_E, g_E\rangle_\mathrm{Pet}\cdot\langle P, P\rangle_\mathrm{NT}}=8.37758040956...=\frac{8\pi}{3}$$ and the equality (\ref{eq:heegnerAIPF}) indeed holds up to an elementary constant $\frac{8\pi}{3}$. 
\end{example}

\subsection{Unitary Shimura varieties}\label{sec:unit-shim-vari}

Now let us come to the setting of unitary Shimura varieties. Let $E/F$ be a CM extension of a totally real number field. Let $\mathbb{V}$ be a totally definite \emph{incoherent hermitian space over $\mathbb{A}_E$} of rank $n$. Here \emph{incoherent} means that $\mathbb{V}$ is not the base change of a global $E/F$-hermitian space, or equivalently the product of the local Hasse invariants of the local hermitian spaces $\mathbb{V}_v:=\mathbb{V} \otimes_{\mathbb{A}_F}F_v$ is $$\prod_{v} \varepsilon(\mathbb{V}_v)=-1.$$ Picking any place $w|\infty$ of $F$ gives a nearby global $E/F$-hermitian space $V$ such that
        \begin{center}
          $V_v\cong\mathbb{V}_v,\text{ if } v\ne w$, \quad{} but $V_{w}$ has
          signature $(n-1,1)$.
        \end{center} Associated to any open compact subgroup $K\subseteq\U(V)(\mathbb{A}_F^\infty)$, we have a \emph{unitary Shimura variety} $X=\Sh_{\U(V)}$ (cf. \cite{Zha19,Gro20}), which has a smooth canonical model of dimension $n-1$ over $E$ (viewed as a subfield of $\mathbb{C} $ via the embedding induced by a place above $w$) and admits complex uniformization $$X(\mathbb{C})=\U(V)(F)\backslash[ \mathbb{D}\times \U(V)(\mathbb{A}_F^\infty)/K],$$  where $$\mathbb{D}:=\{z\in \mathbb{C}^{n-1}: |z|< 1\}\cong \frac{\UU(n-1,1)}{\UU(n-1)\times\UU(1)}$$ is the hermitian symmetric domain associated to $\U(V_\infty)$. 

        We remark that $X$ is a Shimura variety \emph{of abelian type} (rather than of PEL or Hodge type). Unlike Shimura varieties of  PEL type associated to unitary similitude groups, it lacks a good moduli description in terms of abelian varieties with additional structures and thus it is technically more difficult to study. Nevertheless, its \'etale cohomology and $L$-function will be computed in terms of automorphic forms in the forthcoming work of Kisin--Shin--Zhu \cite{Kisina}, under the help of the endoscopic classification for unitary groups due to Mok \cite{Mok15} and Kaletha--Minguez--Shin--White \cite{KMSW14}. In particular, the analytic side of \eqref{eq:nonvanishing} indeed makes sense.

  \subsection{Arithmetic theta lifting} From now on assume that $n=2m$ is even. Let $W=E^{n}$ be the standard $E/F$-skew-hermitian space with matrix $\left(\begin{smallmatrix}  0& 1_m\\ -1_m& 0 \end{smallmatrix}\right)$. Its unitary group $\U(W)$ is a quasi-split unitary group of rank $n$. Let  $\pi$ be a cuspidal automorphic representation of $\U(W)(\mathbb{A}_F)$.

  \begin{assumption}\label{ass:pi} We impose the following (mild) local assumptions on $E/F$ and $\pi$.
    \begin{enumerate}
      \item $E/F$ is split at all 2-adic places and $F\ne \mathbb{Q}$. If $v\nmid \infty$ is ramified in $E$, then $v$ is unramified over $\mathbb{Q}$. Assume that $E/\mathbb{Q}$ is Galois or contains an imaginary quadratic field (for simplicity).
      \item For $v|\infty$, $\pi_v$ is the holomorphic discrete series with Harish-Chandra parameter $\{\frac{n-1}{2}, \frac{n-3}{2},\allowbreak\ldots,\frac{-n+3}{2}, \frac{-n+1}{2}\}$.
      \item For $v\nmid\infty$, $\pi_v$ is tempered.
      \item For $v\nmid\infty$ ramified in $E$, $\pi_v$ is spherical with respect to the stabilizer of $O_{E_v}^{2m}$.
      \item\label{item:inert} For $v\nmid\infty$ inert in $E$, $\pi_v$ is unramified or almost unramified with respect to the stabilizer of $O_{E_v}^{2m}$. If $\pi_v$ is almost unramified, then $v$ is unramified over $\mathbb{Q}$.
      \end{enumerate}
    \end{assumption}

    \begin{remark}
           In Assumption \ref{ass:pi} (\ref{item:inert}), $\pi_v$ is \emph{almost unramified} means that it has a nonzero Iwahori-fixed vector and its Satake parameter contains $\{q_v, q_v^{-1}\}$ and $2m-2$ complex numbers of norm 1. Equivalently, the theta lift of $\pi_v$ to the non-quasi-split unitary group of same rank is spherical with respect to the stabilizer of an almost self-dual lattice (see \cite{Liu21a}).
    \end{remark}

    Let $S_\pi=\{v \text{ inert}: \pi_v\text{ almost unramified}\}$. Then under Assumption \ref{ass:pi}, the global root number for the (complete) standard $L$-function $L(s,\pi)$ (cf. \cite{Yam14}) equals $$\varepsilon(\pi)=(-1)^{|S_\pi|}\cdot (-1)^{m[F:\mathbb{Q}]}$$ by the epsilon dichotomy for unitary groups due to Harris--Kudla--Sweet \cite{Harris1996} and Gan--Ichino \cite[Theorem 11.1]{GI14}. When $\ord_{s=1/2}L(s,\pi)=1$, we have $\varepsilon(\pi)=-1$, and hence there is a canonical choice of totally definite incoherent space $\mathbb{V}=\mathbb{V}_\pi$  of rank $n$ such that for $v\nmid\infty$,
          \begin{center}
            $\varepsilon(\mathbb{V}_v)=-1$ exactly for $v\in S_\pi$.
          \end{center} Let $X$ be the associated unitary Shimura variety  of dimension $n-1=2m-1$ over $E$. The assumption $F\ne \mathbb{Q}$ implies that $X$ is projective.

          \begin{definition}\label{def:arithm-theta-lift}
                      Assuming Kudla's modularity Conjecture \ref{conj:modularity}, then Kudla's generating series $Z_\varphi(\sz)$ of codimension $\i$ special cycles is a hermitian modular form on the hermitian half plane $\mathcal{H}_\i$ valued in $\Ch^\i(X)_\mathbb{C}$. For holomorphic forms $\phi\in \pi$, define the \emph{arithmetic theta lift} by taking the Petersson inner product on $\mathcal{H}_\i$,  $$\Theta_\varphi(\phi):=(Z_\varphi(\sz), \phi)_\text{Pet}\in \Ch^m(X)_\mathbb{C},$$  Moreover, under Assumption \ref{ass:pi}, $\Theta_\varphi(\phi)$ is in fact cohomologically trivial and lies in the $\pi$-isotopic Chow group $\Ch^m(X)_{\pi}^0$ (see \cite[Proposition 6.10]{LL2020}). 
           \end{definition}

\subsection{Arithmetic inner product formula}

Now we are ready to state the arithmetic inner product formula for unitary Shimura varieties.

          \begin{theorem}[with Liu \cite{LL2020, LL2021}]\label{thm:AIPF}
            Let $\pi$ be a cuspidal automorphic representation of $\U(W)(\mathbb{A}_F)$ satisfying Assumption \ref{ass:pi}. Assume that $\varepsilon(\pi)=-1$. Assume that Kudla's modularity Conjecture~\ref{conj:modularity} holds. Then for any holomorphic $\phi\in \pi$ and $\varphi\in \sS(\mathbb{V}(\mathbb{A}_F^\infty)^m)$, the following identity holds (up to simpler factors depending on $\phi$ and $\varphi$),
            \begin{equation}
              \label{eq:AIPF}
              L'(1/2,\pi)\doteq \langle \Theta_{\varphi}(\phi), \Theta_\varphi(\phi)\rangle_\mathrm{BB}.             \end{equation}
          \end{theorem}

          \begin{remark}
            The simpler factors can further be made explicit. For example, if
    \begin{itemize}
    \item $\pi$ is unramified or almost unramified at all finite places,
    \item $\phi\in\pi$ is a holomorphic newform such that $(\phi,\overline{\phi})_\pi=1$,
    \item $\varphi_v$ is the characteristic function of self-dual or almost self-dual lattices at all finite places $v$.
    \end{itemize}
    Then we have
    \begin{equation}\label{eq:explicit}
\frac{L'(1/2,\pi)}{\prod_{i=1}^{2m}L(i, \eta_{E/F}^i)}C_m^{[F:\mathbb{Q}]}\prod_{v\in S_\pi}\frac{q_v^{m-1}(q_v+1)}{(q_v^{2m-1}+1)(q_v^{2m}-1)}=(-1)^m\langle \Theta_\varphi(\phi),\Theta_\varphi(\phi)\rangle_\text{BB},
\end{equation}
where $C_m=2^{-2m}\pi^{m^2}\frac{\Gamma(1)\cdots\Gamma(m)}{\Gamma(m+1)\cdots\Gamma(2m)}$.

Notice that the Grand Riemann Hypothesis predicts that $L'(1/2,\pi)\ge0$, while Beilinson's Hodge index conjecture (\cite[Conjecture 5.5]{Beuilinson1987}) predicts that $(-1)^m\langle \Theta_\varphi(\phi),\Theta_\varphi(\phi)\rangle_\text{BB}\ge0$. It is a good reality check that these two (big) conjectures are compatible with \eqref{eq:explicit}.         \end{remark}

Without assuming Kudla's modularity conjecture, we cannot define $\Theta_\varphi(\phi)$ but we may still obtain unconditional nonvanishing results on Chow groups as predicted by the Beilinson--Bloch conjecture in \eqref{eq:nonvanishing}.

\begin{theorem}[with Liu \cite{LL2020, LL2021}]\label{thm:BB}
  Let $\pi$ be a cuspidal automorphic representation of $\U(W)(\mathbb{A}_F)$ satisfying Assumption \ref{ass:pi}. Let $\Ch^{m}(X)^0_{\mathfrak{m}_\pi}$ the localization of $\Ch^{m}(X)^0_\mathbb{C}$ at the maximal ideal $\mathfrak{m}_\pi$ of the spherical Hecke algebra (away from all ramification) associated to $\pi$. Then the implication
  \begin{equation}
    \label{eq:BBmpi}
    \ord_{s=1/2}L(s,\pi)=1\Longrightarrow \rank\Ch^m(X)^0_{\mathfrak{m}_\pi}\ge1
  \end{equation}
  holds when the level subgroup $K\subseteq \U(V)(\mathbb{A}_F^\infty)$ is sufficiently small.
\end{theorem}

\begin{remark}
  The implication analogous to (\ref{eq:BBmpi}) was previously known for several low dimensional $X$ including:
  \begin{itemize}
  \item   modular curves (Gross--Zagier \cite{Gross1986}),
  \item Shimura curves (Zhang \cite{Zha01}, Kudla--Rapoport--Yang \cite{Kudla2006}, Yuan--Zhang--Zhang \cite{YZZ09}, Liu \cite{Liu2011a}),
  \item  certain $\U(2)\times \U(3)$ Shimura threefolds when  $\pi$ is endoscopic (Xue \cite{Xue19}).
  \end{itemize}
\end{remark}

\subsection{Symmetric power $L$-functions of elliptic curves}

We illustrate Theorem \ref{thm:BB} by an example coming from symmetric power $L$-functions of elliptic curves, which is particularly attractive in view of recent progress on the symmetric power functionality by Newton--Thorne \cite{NT21}. 

\begin{example} Let $A/F$ be a modular elliptic curve without complex multiplication such that
        \begin{itemize}
      \item $\Sym^{2m-1} A$ is modular.
      \item $A$ has bad reduction only at places $v$ split in $E$.
      \end{itemize}
      Assume that $E/F$ satisfies Assumption \ref{ass:pi}. Then there exists $\pi$ satisfying Assumption \ref{ass:pi} such that $$L(s+1/2,\pi)=L(\Sym^{2m-1}A_E,s+m).$$ As $S_\pi=\varnothing$ and $\varepsilon(\pi)=(-1)^{m[F:\mathbb{Q}]}$, Theorem \ref{thm:BB} applies to $\pi$ when $m[F: \mathbb{Q}]$ is odd (e.g., when $m=3$ and $F$ is a totally real cubic field). 
    \end{example}

\subsection{Summary} We end our discussion by the analogy between geometric and arithmetic formulas in Table \ref{tab:summary}.
\begin{center}
  \begin{table}[h]
    \centering
    \begin{tabular}{c|c}
      Automorphic/Geometric & Arithmetic\\\hline
      Hurwitz formula (\ref{eq:hurwitz}) & Gross--Keating formula (\ref{eq:grosskeating})\\\hline
      Geometric Siegel--Weil formula (\ref{eq:kudla}) & Arithmetic Siegel--Weil formula (\ref{eq:ASW}) \\  (Kudla's formula) & (Kudla--Rapoport Conjecture)\\
      \hline
      Rallis inner product formula (\ref{eq:rallis}) & Arithmetic inner product formula  (\ref{eq:AIPF})\\ & (Gross--Zagier formula in higher dimension)
    \end{tabular}
    \caption{Summary}\label{tab:summary}
  \end{table}
  \end{center}

\bibliographystyle{alpha}
\bibliography{KR}

\end{document}